\PassOptionsToPackage{unicode}{hyperref}
\PassOptionsToPackage{hyphens}{url}
\PassOptionsToPackage{dvipsnames,svgnames,x11names}{xcolor}
\documentclass[
  11pt,
]{article}
\usepackage{amsmath,amssymb}
\usepackage{iftex}
\ifPDFTeX
  \usepackage[T1]{fontenc}
  \usepackage[utf8]{inputenc}
  \usepackage{textcomp} 
\else 
  \usepackage{unicode-math} 
  \defaultfontfeatures{Scale=MatchLowercase}
  \defaultfontfeatures[\rmfamily]{Ligatures=TeX,Scale=1}
\fi
\usepackage{lmodern}
\ifPDFTeX\else
\fi
\IfFileExists{upquote.sty}{\usepackage{upquote}}{}
\IfFileExists{microtype.sty}{
  \usepackage[]{microtype}
  \UseMicrotypeSet[protrusion]{basicmath} 
}{}
\makeatletter
\@ifundefined{KOMAClassName}{
  \IfFileExists{parskip.sty}{%
    \usepackage{parskip}
  }{
    \setlength{\parindent}{0pt}
    \setlength{\parskip}{6pt plus 2pt minus 1pt}}
}{
  \KOMAoptions{parskip=half}}
\makeatother
\usepackage{xcolor}
\usepackage[margin=1in]{geometry}
\usepackage{longtable,booktabs,array}
\usepackage{calc} 
\usepackage{etoolbox}
\makeatletter
\patchcmd\longtable{\par}{\if@noskipsec\mbox{}\fi\par}{}{}
\makeatother
\IfFileExists{footnotehyper.sty}{\usepackage{footnotehyper}}{\usepackage{footnote}}
\makesavenoteenv{longtable}
\usepackage{graphicx}
\makeatletter
\def\maxwidth{\ifdim\Gin@nat@width>\linewidth\linewidth\else\Gin@nat@width\fi}
\def\maxheight{\ifdim\Gin@nat@height>\textheight\textheight\else\Gin@nat@height\fi}
\makeatother
\setkeys{Gin}{width=\maxwidth,height=\maxheight,keepaspectratio}
\makeatletter
\def\fps@figure{htbp}
\makeatother
\providecommand{\tightlist}{%
  \setlength{\itemsep}{0pt}\setlength{\parskip}{0pt}}
\newlength{\cslhangindent}
\newlength{\csllabelwidth}
\newlength{\cslentryspacingunit} 
\newenvironment{CSLReferences}[2] 
 {
  \setlength{\parindent}{0pt}
  \ifodd #1
  \let\oldpar\par
  \def\par{\hangindent=\cslhangindent\oldpar}
  \fi
  \setlength{\parskip}{#2\cslentryspacingunit}
 }%
 {}
\usepackage{calc}

\ifLuaTeX
  \usepackage{selnolig}  
\fi
\IfFileExists{bookmark.sty}{\usepackage{bookmark}}{\usepackage{hyperref}}
\IfFileExists{xurl.sty}{\usepackage{xurl}}{} 
\hypersetup{
  pdftitle={Exact Minimum Field Partition for One-Step Prediction in a Finite Five-Agent ABCW Model},
  pdfauthor={Takashi Inoue Shizusawa},
  pdfkeywords={agent-based model, predictive partition, state
aggregation, graph coloring, coarse-graining, exact prediction},
  colorlinks=true,
  linkcolor={blue},
  filecolor={Maroon},
  citecolor={blue},
  urlcolor={blue},
  pdfcreator={LaTeX via pandoc}}

\title{Exact Minimum Field Partition for One-Step Prediction in a Finite
Five-Agent ABCW Model}
\author{Takashi Inoue Shizusawa\\
\small Independent Researcher, Japan}
\date{}

\begin{document}
\maketitle
\begin{abstract}
Which distinctions in a present state must be retained to predict a
specified future exactly? We study this question on a finite reachable
set of a deterministic five-agent ABCW model in which binary agent
actions and strategies coevolve with a weighted directed influence
field. The dataset is generated from four initial-field families and all
4,096 action-strategy initial conditions, yielding 56,536 transitions
and 2,562 distinct current fields. We retain the current action and seek
the coarsest field-only partition that uniquely determines the complete
one-step-ahead field on every observed transition. Nine natural feature
families and all 511 of their nonempty combinations provide strong
approximate predictors but no exact nontrivial field compression. We
therefore construct an incompatibility graph whose vertices are current
fields and whose edges join fields that produce different next fields
under a shared observed action. A proper 692-coloring gives a
constructive upper bound. An independent anchor-action decomposition,
followed by exact solution of the 50 nontrivial induced subgraphs, gives
the matching lower bound. Hence the minimum number of field classes is
exactly 692. This number applies only to the specified finite reachable
set, field-only compression, retained current actions, one-step target,
and exact-prediction requirement; it is not a universal macrostate count
for ABCW systems.
\end{abstract}


\vspace{0.75em}

\hypertarget{introduction}{%
\section{Introduction}\label{introduction}}

In complex dynamical systems, a more detailed description of the present
does not necessarily provide a more useful description for predicting
the future. A complete microscopic state may be sufficient to determine
future evolution, but not every distinction encoded in that state need
be relevant to a specified predictive task.

The central question of this study is therefore:

\begin{quote}
\textbf{Which distinctions in the present must be retained in order to
predict the future?}
\end{quote}

Equivalently:

\begin{quote}
\textbf{How much of the present can be forgotten without losing the
specified future?}
\end{quote}

This question is related to several established lines of research,
including coarse-graining and state aggregation, lumpability in Markov
chains, model minimization through bisimulation, causal states in
computational mechanics, and state minimization for incompletely
specified finite-state machines. These approaches share a broad
objective: to reduce a state description while preserving the
distinctions required for dynamics or prediction. They differ, however,
in what must be preserved and in the conditions imposed on the reduced
representation.

The present study does not propose a new general theory encompassing
these frameworks. Instead, it considers a finite and explicitly
computable agent-based dynamical system and determines the minimum
number of distinctions required to preserve a specified one-step
prediction exactly.

The system studied here is a five-agent ABCW model. Each agent has a
binary action and a binary strategy, while the influence relations among
agents are represented by a time-dependent weight matrix \({W}\).
Actions are updated by reference to local influence relations, and game
outcomes feed back into both strategy adaptation and field evolution.

The payoff rule has a Minority-Game-type competitive structure. The ABCW
model is not intended to reproduce the standard Minority Game itself.
Rather, the rule is used to introduce an environment in which all agents
cannot simultaneously succeed in the same sense. The relationship of
this construction to the Minority Game and the El Farol problem is
discussed in Section 2.

ABCW is therefore a finite agent system in which agent behavior and the
interaction field coevolve. We represent the field by its deviation from
the fully symmetric reference field \({W^\ast}\):

\[{\Delta W=W-W^\ast.}\]

For the five-agent model, we use four initial fields and exhaust all
action and strategy initial conditions for each field, giving a total of
4,096 initial conditions. Each trajectory is followed until the complete
state recurs. The resulting finite reachable dataset contains 56,536
transitions, 2,562 distinct current fields \({\Delta W}\), and 11,202
distinct \({(a,\Delta W)}\) states.

Within this finite system, an observation is treated not merely as the
selection of a variable but as a mapping

\[{P:\Omega\longrightarrow Y,}\]

which specifies which distinctions among complete states are retained
and which are discarded.

A key distinction is that closure of the observed dynamics and
prediction of a specified target are different requirements. For an
observation \({P}\), requiring

\[{P(X_t)\longrightarrow P(X_{t+1})}\]

to be single-valued is a self-closure condition. By contrast, for a
specified prediction target \({Z}\), requiring

\[{P(X_t)\longrightarrow Z}\]

to be single-valued is a target-prediction condition. The primary
objective of this study is the latter.

Specifically, we retain the current action \({a_t}\) and replace the
current field \({\Delta W_t}\) with a coarser field-class label. We then
ask whether this compressed description still determines the complete
one-step-ahead field \({\Delta W_{t+1}}\) uniquely.

Because the object of compression is the field description, the action
\({a_t}\) is retained as an uncompressed conditioning variable. We do
not simultaneously optimize a compression of the action space. This
field-only restriction is part of the scope of the minimization problem
solved in this paper.

Let

\[{B:\Delta W\longrightarrow\mathcal C}\]

be a field-only mapping, and define

\[{P_B(X_t)=\bigl(a_t,B(\Delta W_t)\bigr).}\]

We require

\[{\bigl(a_t,B(\Delta W_t)\bigr)\longrightarrow\Delta W_{t+1}}\]

to be single-valued over all observed transitions. We call the ability
of an observation to determine the specified target without ambiguity
\textbf{predictive sufficiency}.

We first examine field compression based on human-designed network
features, including norms, local structural descriptors, and
out-strength. Some of these features achieve high predictive
performance. Within the candidate family tested, however, no feature
combination both merges genuinely distinct values of \({\Delta W}\) and
predicts the complete one-step-ahead field exactly.

This failure does not imply that exact prediction requires all 2,562
observed fields to remain distinct. It shows only that the tested
feature family does not identify an exact nontrivial compression. We
therefore remove the restriction that the partition must be specified in
advance through human-designed features and instead derive the necessary
distinctions directly from the observed dynamics.

The finite dataset can be read as a partially specified input-output
relation

\[{(\Delta W,a)\longmapsto\Delta W',}\]

where \({\Delta W}\) is the current field, \({a}\) is the action
condition, and \({\Delta W'}\) is the complete one-step-ahead field. If
two current fields produce different next fields under a common action
condition, they cannot belong to the same field class.

Accordingly, we construct an incompatibility graph \({G}\) whose 2,562
vertices are the distinct observed current fields. Two fields are
connected by an edge if, under at least one action condition observed
for both fields, they produce different next fields. A field-only
partition that preserves exact one-step prediction is then equivalent to
a proper coloring of \({G}\).

The reduced field-class label is not required to serve recursively as
the internal state at the next time step. The target being preserved is
the complete next field \({\Delta W_{t+1}}\), rather than the successor
class label. Consequently, the recursive successor-closure condition
encountered in classical minimization of incompletely specified
finite-state machines is not imposed here. For this depth-1 prediction
problem, the minimum number of field classes is therefore given directly
by

\[{|\operatorname{Im}B_{\min}|=\chi(G).}\]

Compatibility need not itself be assumed to form an equivalence
relation. Pairwise distinctions that cannot be merged are encoded
directly as incompatibility edges, and a proper coloring separates every
incompatible pair. The relationship between this formulation and
classical state-minimization problems is developed in Section 2.

The main result is

\[{\chi(G)=692.}\]

Thus, for the finite reachable set analyzed here and for the specified
task of exact one-step field prediction, the 2,562 distinct current
fields can be compressed to 692 field classes:

\[{2562\longrightarrow692,}\]

while preserving exact prediction of the complete next field.

Figure 1 summarizes the predictive-partition problem and the resulting
reduction.

\begin{figure}
\centering
\includegraphics[trim=0 0 0 24bp,clip]{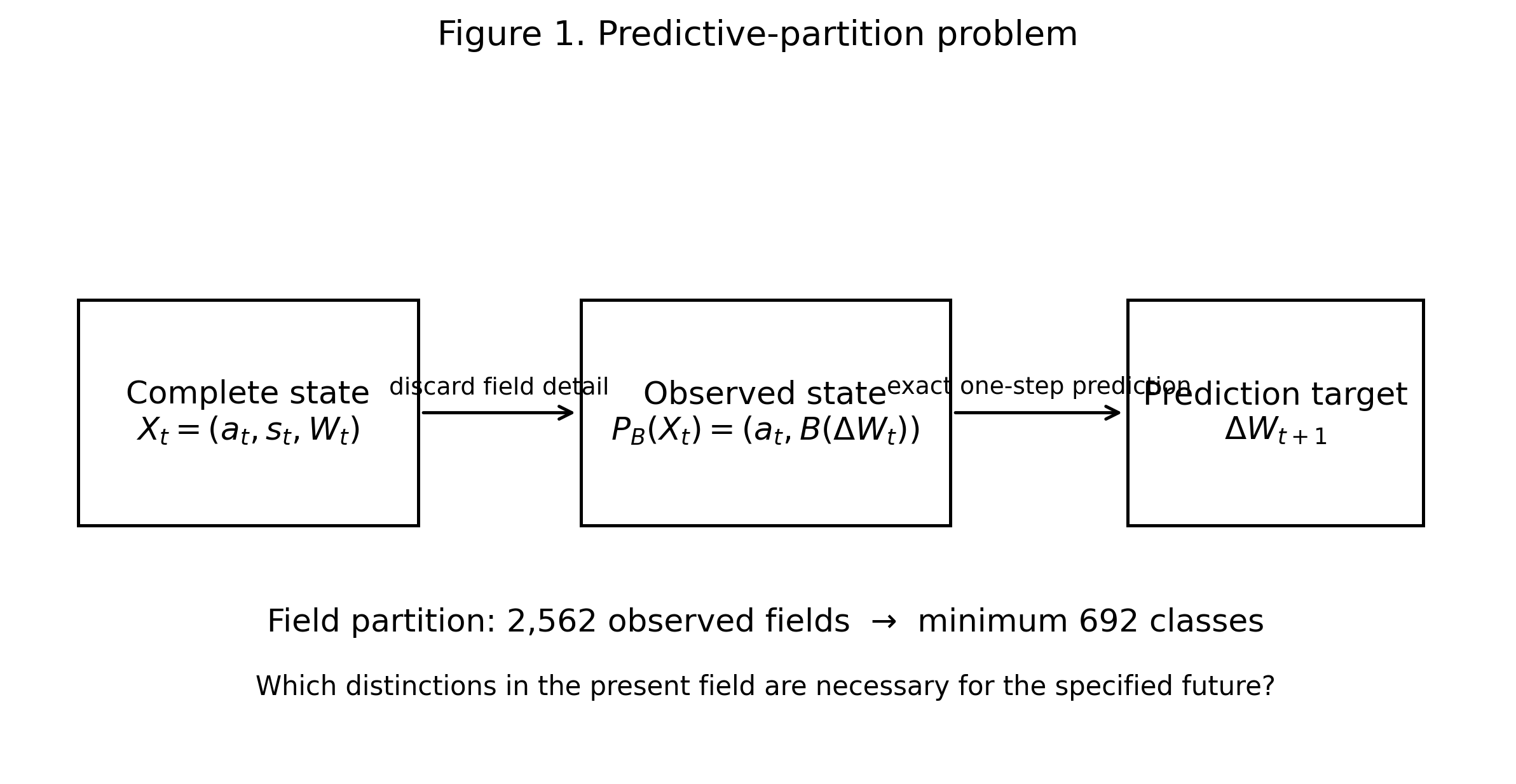}
\caption{Predictive-partition problem. The complete ABCW state
is compressed by retaining the action and replacing the full field
deviation with a field-only partition label. The prediction target is
the exact one-step-ahead field deviation. The main result reduces 2,562
observed fields to a minimum of 692 classes.}
\end{figure}

Measured by the number of field classes, the corresponding compression
rate is

\[{1-\frac{692}{2562}\simeq72.99\%.}\]

This is a reduction in the number of field classes; it should not be
interpreted as a 72.99\% reduction in information content.

The value 692 is not merely the output of a coloring heuristic. A lower
bound is obtained by fixing the anchor action

\[{a=(-,-,-,-,-).}\]

Under this action, the 1,239 observed current fields produce 625
distinct next fields, so at least 625 field classes are required. Within
50 of these 625 anchor groups, additional incompatibilities appear when
other action conditions are taken into account. The chromatic numbers of
the corresponding induced subgraphs are solved exactly by combining
clique lower bounds with colorability searches. These refinements
require 67 additional classes, giving the lower bound

\[{625+67=692.}\]

For the matching upper bound, we explicitly construct a proper
692-coloring of the incompatibility graph. Mapping this coloring back to
all 56,536 observed transitions verifies that every observed pair
\({(a_t,B(\Delta W_t))}\) determines a unique \({\Delta W_{t+1}}\). The
upper and lower bounds therefore coincide, establishing
\({\chi(G)=692}\) exactly.

This result does not imply that 692 is a universal number of macrostates
for ABCW systems. It depends on the five-agent model, the four
initial-field families---Baseline, Hub, Local, and Hub+Local---the
finite reachable set generated from them, the retention of action as an
uncompressed conditioning variable, the field-only restriction, and the
requirement of exact prediction of the complete one-step-ahead field. We
therefore do not claim that 692 is the minimum state count of ABCW in
general.

The broader point is that microscopic distinctness and predictive
relevance need not coincide. In this finite ABCW system, fields that are
microscopically different can sometimes be assigned to the same
predictive class without losing the specified future. At the same time,
simple network features and geometric proximity between fields do not
fully identify which distinctions must be preserved.

The resulting partition occupies an intermediate descriptive level
between a complete microscopic representation and a small set of simple
macroscopic statistics: it retains exactly the distinctions required by
the specified prediction task on the observed finite system. This does
not by itself make the 692 classes an interpretable set of macroscopic
variables, but it provides an exact target against which candidate
coarse descriptions can be evaluated.

The contributions of this study are threefold.

First, for a finite ABCW system, we distinguish self-closure from target
prediction and formulate the question of which present distinctions must
be retained as an explicit observation-minimization problem.

Second, we distinguish compression by preselected network features from
a predictive partition derived from future responses, and we formulate
the resulting field-only, action-conditioned, depth-1 problem as a
graph-coloring problem.

Third, for the finite reachable set considered here, we establish that
the exact minimum number of field classes is 692 by providing both a
constructive upper bound and an independent matching lower bound. The
result demonstrates that not all 2,562 microscopic field distinctions
are required for exact one-step field prediction.

These contributions do not assert novelty for the general theories of
state minimization, lumpability, bisimulation, or causal states. Rather,
the study formulates a specifically scoped prediction problem for a
finite coevolving agent--field system and solves its exact minimum on a
fully inspectable dataset.

The remainder of the paper is organized as follows. Section 2 positions
the study in relation to the El Farol problem, the Minority Game, state
aggregation, lumpability, bisimulation, computational mechanics, and
minimization of incompletely specified finite-state machines. Section 3
defines the ABCW model and the finite reachable dataset. Section 4
formalizes observation, self-closure, predictive sufficiency, and the
prediction target used in this paper. Section 5 evaluates compression
based on natural network features. Section 6 constructs the
incompatibility graph and determines the exact minimum field-only
partition. Section 7 analyzes the internal structure of the 692 classes.
Section 8 discusses the interpretation, scope, and limitations of the
result. Section 9 concludes the paper.

\hypertarget{related-work-and-positioning}{%
\section{Related Work and
Positioning}\label{related-work-and-positioning}}

\hypertarget{competitive-agent-models-el-farol-and-the-minority-game}{%
\subsection{Competitive agent models: El Farol and the Minority
Game}\label{competitive-agent-models-el-farol-and-the-minority-game}}

The competitive payoff structure of the ABCW model shares part of its
motivation with the El Farol Bar problem and the Minority Game.

The El Farol Bar problem introduced by Arthur
(\protect\hyperlink{ref-arthur1994}{1994}) considers a setting in which
several agents seek access to the same limited opportunity or resource,
so that the desirability of an action depends on the choices of other
agents. Attendance is advantageous when relatively few agents attend,
but its value declines when too many make the same choice. Consequently,
the rationality of an action cannot be fixed independently of the
behavior of the population.

The Minority Game introduced by Challet and Zhang
(\protect\hyperlink{ref-challet1997}{1997}) formulates this competitive
structure as a simplified binary-choice game. At each time step, agents
select one of two alternatives, and agents in the minority receive a
payoff. The value of an action is therefore not intrinsic to that action
but depends on how many other agents choose it.

ABCW employs this Minority-Game-type competition. Each of five agents
takes a binary action, with agents in the minority designated as winners
and those in the majority as losers. This creates a competitive payoff
environment in which all agents cannot simultaneously succeed in the
same sense.

ABCW is not intended, however, to reproduce or extend the standard
Minority Game, nor is it designed to analyze its characteristic
collective phenomena. Standard Minority Game research commonly focuses
on agents' strategy sets, information histories, adaptation, collective
efficiency, and aggregate fluctuations. By contrast, the central object
in the present study is the time-dependent field \({W}\) representing
influence relations among agents.

In ABCW, game outcomes feed back not only into agent scores but also
into the evolution of the field. The Minority-Game-type rule is
therefore a component that generates the coupled process

\[{\text{action}\longrightarrow\text{outcome}\longrightarrow\text{field update}\longrightarrow\text{next action},}\]

rather than the final object of analysis. The primary question is not
which strategy dominates in a Minority Game or how efficiently the
population coordinates. Instead, for the field dynamics generated by
this competitive agent system, we ask:

\begin{quote}
\textbf{Which distinctions among current fields must be retained to
predict the complete one-step-ahead field exactly?}
\end{quote}

The El Farol problem and the Minority Game thus provide the
agent-modeling background for the competitive dynamics of ABCW. They are
not the direct theoretical basis of the predictive state-reduction
problem addressed here. The more immediate connections are to state
aggregation and finite-state minimization.

\hypertarget{state-aggregation-and-the-question-of-relevant-distinctions}{%
\subsection{State aggregation and the question of relevant
distinctions}\label{state-aggregation-and-the-question-of-relevant-distinctions}}

The central question of this study---which distinctions in the present
must be retained for prediction---belongs to a broad class of problems
concerned with aggregating a state space while preserving properties
required for dynamics or prediction.

A classical example is lumpability in Markov chains. In the formulation
of Kemeny and Snell (\protect\hyperlink{ref-kemeny1960}{1960}), a
partition of the state space into lumps yields a closed Markov process
on the aggregated states when, for any two states in the same lump, the
transition probabilities into every destination lump coincide.
Lumpability therefore provides conditions under which microscopic states
can be merged while retaining well-defined reduced dynamics.

The self-closure condition introduced in Section 4 has a related
motivation in deterministic finite dynamics. If

\[{P(X)=P(X')}\]

implies

\[{P(\Phi(X))=P(\Phi(X')),}\]

then two complete states assigned to the same observed state also evolve
to the same observed state at the next time step. Here \({X}\) is a
complete state, \({\Phi}\) is the one-step update map on complete
states, and \({P:\Omega\to Y}\) is a general observation map on the
complete state space. This general map \({P}\) should be distinguished
from the field-only map \({B:\Delta W\to\mathcal C}\) that is later
minimized.

Self-closure is not, however, the primary requirement in this paper. We
distinguish prediction of the next observed state,

\[{P(X_{t+1}),}\]

from prediction of a separately specified target,

\[{Z=g(X_{t+1}).}\]

The required condition is therefore not necessarily closure of

\[{P(X_t)\longrightarrow P(X_{t+1}),}\]

but single-valuedness of

\[{P(X_t)\longrightarrow Z.}\]

This distinction between closed reduced dynamics and sufficient
information for a specified target is the starting point of our
state-aggregation problem.

\hypertarget{bisimulation-and-model-minimization}{%
\subsection{Bisimulation and model
minimization}\label{bisimulation-and-model-minimization}}

The idea of merging states while preserving future behavior is also
closely related to bisimulation and bisimulation-based model
minimization.

Givan, Dean, and Greig (\protect\hyperlink{ref-givan2003}{2003}), for
example, study equivalence notions for Markov decision processes and use
bisimulation-based aggregation to construct reduced MDPs that preserve
optimal policies of the original process. The common principle is that
states are grouped not because they are superficially similar, but
because their distinction is unnecessary for the future behavior that
must be preserved.

The system and preservation target considered here are different. ABCW
is not an MDP; on the stored finite reachable set, it has deterministic
update rules. Nor do we minimize the complete ABCW state. Writing the
complete state as

\[{X=(a,s,W),}\]

we retain the action \({a}\) in the observation and compress only the
field deviation \({\Delta W}\) through

\[{B:\Delta W\longrightarrow\mathcal C.}\]

Moreover, the preserved target is neither an optimal policy nor a reward
structure nor the full reduced dynamics. It is the complete
one-step-ahead field \({\Delta W_{t+1}}\). We therefore do not describe
the resulting minimum partition as a bisimulation quotient.

\hypertarget{computational-mechanics-and-causal-states}{%
\subsection{Computational mechanics and causal
states}\label{computational-mechanics-and-causal-states}}

Computational mechanics provides another closely related predictive
perspective. Shalizi and Crutchfield
(\protect\hyperlink{ref-shalizi2001}{2001}) define causal states by
identifying past histories that induce the same conditional probability
distribution over futures. The resulting \({\epsilon}\)-machine is
characterized as a minimal representation retaining the information
required for prediction.

The conceptual connection to the present study is clear: distinctions
among complete histories or microscopic states are retained only when
they matter for the future. Similarly, the 692 classes derived here are
defined not by static similarity in network structure but by the
one-step field responses generated under action conditions.

The two constructions should nevertheless not be identified. Causal
states are, in principle, equivalence classes of histories defined by
conditional distributions over future sequences. The present study
instead partitions current fields \({\Delta W}\) in a finite reachable
set, with the prediction target restricted to the complete
one-step-ahead field \({\Delta W_{t+1}}\).

The ABCW data are also partially specified: not every action condition
is observed for every current field. As a result, the relation ``no
conflict occurs under the observed common conditions'' need not be a
transitive equivalence relation. We therefore do not call the 692
classes causal states or an \({\epsilon}\)-machine. The relevant
connection is the shared objective of preserving future-relevant
distinctions.

\hypertarget{incompletely-specified-finite-state-machine-minimization}{%
\subsection{Incompletely specified finite-state machine
minimization}\label{incompletely-specified-finite-state-machine-minimization}}

The most direct mathematical correspondence is with minimization of
incompletely specified finite-state machines. A classical treatment is
given by Paull and Unger (\protect\hyperlink{ref-paull1959}{1959}). When
outputs or transitions are not specified for every state-input pair,
state equivalence for a fully specified machine cannot be applied
directly. Minimization methods instead examine whether states are
mutually compatible and construct compatible classes that permit a
reduction in the number of states. This remains a standard
state-reduction problem in finite-state machine theory
(\protect\hyperlink{ref-paull1959}{Paull and Unger 1959};
\protect\hyperlink{ref-kohavi2009}{Kohavi and Jha 2009}).

In classical minimization of an incompletely specified machine,
compatibility of current input-output behavior is generally not
sufficient. If the reduced machine is to operate recursively, successor
states must also be consistent with the selected reduced classes. The
minimization problem therefore involves a cover satisfying both
compatibility and closure conditions.

The ABCW dataset has a related structure when viewed as partially
specified input-output behavior. Taking \({\Delta W}\) as the current
field, \({a}\) as an input condition, and \({\Delta W'}\) as the
complete next-field output, the observations define a partial table

\[{(\Delta W,a)\longmapsto\Delta W'.}\]

Suppose that two current fields \({\Delta W_i}\) and \({\Delta W_j}\)
are both observed under the same action condition \({a}\) but generate
different next fields,

\[{\Delta W_i'\neq\Delta W_j'.}\]

The two current fields cannot then be assigned to the same field class.
We represent this constraint by an incompatibility graph whose 2,562
vertices are the distinct observed current fields and whose edges join
pairs that produce different next fields under at least one common
observed action condition.

Adjacent vertices cannot share a field class. A field-only partition
preserving exact one-step prediction therefore gives a proper coloring
of the incompatibility graph. Conversely, in a proper coloring, two
fields assigned the same color have no observed common action condition
under which they produce conflicting next fields. Treating each color as
a field class consequently preserves the single-valued relation

\[{\bigl(a_t,B(\Delta W_t)\bigr)\longrightarrow\Delta W_{t+1}.}\]

The minimum number of classes is therefore exactly

\[{|\operatorname{Im}B_{\min}|=\chi(G).}\]

An important difference from classical incompletely specified FSM
minimization is that the field-class label is not required to function
recursively as the internal state at the next time step. The object
preserved here is the complete next field \({\Delta W_{t+1}}\), not a
successor class label. The requirement is only that the depth-1
target-prediction relation

\[{\bigl(a_t,B(\Delta W_t)\bigr)\longrightarrow\Delta W_{t+1}}\]

be single-valued. No recursive closure condition on successor classes is
therefore required. For the problem defined in this paper, minimization
reduces directly to proper coloring rather than to a closed-cover
construction.

Nor must compatibility itself be assumed to be an equivalence relation.
In partially specified data, the absence of an observed conflict need
not be transitive. We instead encode every pair that cannot be merged as
an incompatibility edge and enforce the required pairwise constraints
through proper coloring.

This simplification is not presented as a new general minimization
theory for incompletely specified finite-state machines. It is a
consequence of the specific preservation target adopted here: because
the ABCW problem is restricted to an external one-step target, the
classical recursive closure requirement is unnecessary, and the minimum
partition is expressible as a pure graph-coloring problem.

\hypertarget{feature-based-coarse-graining-and-predictive-partitioning}{%
\subsection{Feature-based coarse-graining and predictive
partitioning}\label{feature-based-coarse-graining-and-predictive-partitioning}}

We further distinguish two approaches to state compression.

The first computes a set of preselected features from the state:

\[{\Delta W\longrightarrow B_{\mathrm{feature}}(\Delta W).}\]

Examples considered here include norms, local structural descriptors,
hub structure, and out-strength. The advantage of this approach is
interpretability: the resulting observables have meanings chosen in
advance. A low-dimensional feature representation, however, does not
necessarily constitute genuine coarse-graining. On a finite state set,
even a low-dimensional feature vector may uniquely identify every
microscopic state. When features do merge states, they may also discard
distinctions required for prediction.

The second approach does not choose the features in advance. Instead, it
derives the required distinctions from future responses, conceptually
through the response structure

\[{\Delta W\longrightarrow\left(a\longmapsto\Delta W'\right).}\]

The 692-class partition belongs to this second approach. Its classes are
not guaranteed to admit an immediate interpretation in terms of simple
network statistics. In return, the partition is determined directly by
the requirement that the specified predictive target be preserved.

We refer to this property as \textbf{predictive sufficiency}. An
observation \({P}\) is predictively sufficient for a target \({Z}\) on
the dataset under consideration if

\[{P(X_t)\longrightarrow Z}\]

is single-valued, so that \({Z}\) is determined without ambiguity by
\({P(X_t)}\). Predictive sufficiency does not mean merely that a state
can be expressed with a small number of variables. The question is
instead how many distinctions can be discarded while retaining all
distinctions required for the specified target.

\hypertarget{position-of-the-present-study}{%
\subsection{Position of the present
study}\label{position-of-the-present-study}}

The present study has two distinct connections to previous work.

The first concerns the competitive agent dynamics of ABCW. The model
shares with the El Farol problem and the Minority Game a payoff
structure in which an agent's outcome depends on the choices of other
agents and all agents cannot simultaneously succeed in the same sense.
We do not analyze standard Minority Game phenomena or strategy
adaptation as ends in themselves. The competitive rule is used to
generate coupled evolution between agent behavior and the field \({W}\).

The second connection concerns the theoretical positioning of the
predictive state-reduction problem. With lumpability, we share the
concern of obtaining a well-defined description after aggregation, but
we distinguish self-closure from prediction of a separately specified
target. With bisimulation and MDP minimization, we share the principle
of aggregation according to future behavior, but we do not preserve
optimal policies or a complete reduced MDP. With computational
mechanics, we share the objective of retaining future-relevant
distinctions, but the 692 classes are not history-based causal states:
they preserve one-step, action-conditioned field responses on a finite
reachable set. With incompletely specified FSM minimization, the
correspondence is more direct because both problems involve
compatibility and incompatibility in partially specified input-output
behavior; our target restriction makes recursive closure unnecessary.

The most limited characterization of the present study is therefore:

\begin{quote}
\textbf{It is a concrete state-minimization problem that determines the
exact minimum field-only partition preserving a specified one-step
prediction for partially specified input-output data generated by a
finite competitive agent system.}
\end{quote}

We do not locate the novelty of the study in the general ideas of the
Minority Game, state minimization, or predictive representation
themselves. The specific contribution is to formulate the depth-1
prediction problem explicitly for the observed finite reachable set of
the five-agent ABCW model and to establish, for its incompatibility
graph, that

\[{\chi(G)=692}\]

by a constructive upper bound and an independent matching lower bound.
Consequently, the complete one-step-ahead field can be predicted exactly
without retaining all 2,562 microscopically distinct current fields,
whereas no partition into 691 or fewer field classes can preserve that
prediction.

This finite and fully inspectable example makes explicit the difference
between being distinct as a state and needing to remain distinct for a
specified prediction.

\begin{center}\rule{0.5\linewidth}{0.5pt}\end{center}

\hypertarget{abcw-model-and-dataset}{%
\section{ABCW Model and Dataset}\label{abcw-model-and-dataset}}

This section defines the five-agent ABCW model analyzed in this paper
and the finite reachable dataset used in the predictive-partition
problem. The purpose is not to reconstruct the exploratory development
of ABCW, but to fix the state variables, update rules, field
representation, and scope of the subsequent analysis.

\hypertarget{state-variables-and-interaction-network}{%
\subsection{State Variables and Interaction
Network}\label{state-variables-and-interaction-network}}

Let the set of agents be

\[{V=\lbrace1,\ldots,n\rbrace,}\]

with

\[{n=5}\]

throughout this paper. At time \({t}\), each agent \({i}\) has a binary
action

\[{a_i(t)\in\lbrace-1,+1\rbrace}\]

and a binary strategy

\[{s_i(t)\in\lbrace-1,+1\rbrace.}\]

The value \({s_i(t)=+1}\) denotes a trend-following strategy that
follows the reference signal, whereas \({s_i(t)=-1}\) denotes a
contrarian strategy that acts against it.

Influence relations among agents are represented by a fixed topology

\[{E=(E_{ij}),\qquad E_{ij}\in\lbrace0,1\rbrace}\]

and a time-dependent nonnegative weight matrix

\[{W(t)=(w_{ij}(t)).}\]

The condition \({E_{ij}=1}\) indicates that the directed edge
\({i\to j}\) is permitted by the topology, and \({w_{ij}(t)}\)
represents the strength of the influence exerted by agent \({i}\) on
agent \({j}\) at time \({t}\). If its weight reaches zero, the
topological edge itself is not removed as long as \({E_{ij}=1}\).

Because strategy is included as a state variable, the complete state is
defined as

\[{X(t)=(a(t),s(t),W(t)).}\]

Table 1 summarizes the principal variables and their roles.

\hypertarget{table-1.-abcw-variables-and-update-rule-summary}{%
\subsection{Table 1. ABCW variables and update-rule
summary}\label{table-1.-abcw-variables-and-update-rule-summary}}

\begin{longtable}[]{@{}
  >{\raggedright\arraybackslash}p{(\columnwidth - 6\tabcolsep) * \real{0.2500}}
  >{\raggedright\arraybackslash}p{(\columnwidth - 6\tabcolsep) * \real{0.2500}}
  >{\raggedright\arraybackslash}p{(\columnwidth - 6\tabcolsep) * \real{0.2500}}
  >{\raggedright\arraybackslash}p{(\columnwidth - 6\tabcolsep) * \real{0.2500}}@{}}
\toprule\noalign{}
\begin{minipage}[b]{\linewidth}\raggedright
Symbol
\end{minipage} & \begin{minipage}[b]{\linewidth}\raggedright
Meaning
\end{minipage} & \begin{minipage}[b]{\linewidth}\raggedright
Domain / type
\end{minipage} & \begin{minipage}[b]{\linewidth}\raggedright
Role / update
\end{minipage} \\
\midrule\noalign{}
\endhead
\bottomrule\noalign{}
\endlastfoot
\(a_i(t)\) & Action & \(\{-1,+1\}\) & Current binary action \\
\(s_i(t)\) & Strategy & \(\{-1,+1\}\) & Trend-following / contrarian;
loser flips strategy \\
\(u_i(t)\) & Payoff & \(\{-1,+1\}\) for \(n=5\) & Minority \(+1\),
majority \(-1\); all-same gives \(-1\) \\
\(W(t)\) & Influence field & Nonnegative weighted directed matrix &
Outgoing edges updated by source payoff \\
\(\sigma_i(t)\) & Reference signal & \([-1,1]\) & Mean action among
strongest effective incoming references \\
\(\Delta W(t)\) & Field deviation & Matrix & \(W(t)-W^\ast\) \\
\end{longtable}

\hypertarget{payoff-reference-signal-and-update-rules}{%
\subsection{Payoff, Reference Signal, and Update
Rules}\label{payoff-reference-signal-and-update-rules}}

\hypertarget{minority-game-payoff}{%
\subsubsection{Minority-game payoff}\label{minority-game-payoff}}

At each time step, every agent chooses either \({-1}\) or \({+1}\).
Agents in the minority are designated winners and those in the majority
losers. The payoff is

\[{u_i(t)\in\lbrace-1,+1\rbrace,}\]

with \({u_i(t)=+1}\) for a minority agent and \({u_i(t)=-1}\) for a
majority agent. Because \({n=5}\) is odd, a tie between the two binary
actions is impossible. If all agents choose the same action, all are
treated as losers and receive \({u_i(t)=-1}\).

This payoff rule is not intended to reproduce the full standard Minority
Game. It is used to introduce into ABCW a competitive condition under
which all agents cannot win simultaneously.

\hypertarget{reference-signal}{%
\subsubsection{Reference signal}\label{reference-signal}}

The set of agents that agent \({i}\) can reference is

\[{\mathcal N_i=\lbrace j\in V\mid E_{ji}=1\rbrace.}\]

Define the maximum incoming weight available to agent \({i}\) as

\({m_i(t)= \begin{cases} \max_{j\in\mathcal N_i}w_{ji}(t),&\mathcal N_i\neq\varnothing,\cr 0,&\mathcal N_i=\varnothing \end{cases}.}\)

When \({m_i(t)>0}\), let

\[{M_i(t)=\lbrace j\in\mathcal N_i\mid w_{ji}(t)=m_i(t)\rbrace}\]

be the set of agents attaining that maximum, and define the reference
signal as

\[{\sigma_i(t)=\frac{1}{|M_i(t)|}\sum_{j\in M_i(t)}a_j(t).}\]

If several agents share the maximum incoming weight, no single agent is
selected arbitrarily; all tied agents are referenced with equal weight.
If their actions balance, then \({\sigma_i(t)=0}\). We also set
\({\sigma_i(t)=0}\) when no referenceable agent exists or when the
maximum available influence is zero, treating both cases as the absence
of an effective reference signal.

The game outcome \({u_i(t)}\) is thus determined by the actions of all
agents, whereas the reference information \({\sigma_i(t)}\) used for the
next action is determined by local influence relations.

\hypertarget{strategy-adaptation}{%
\subsubsection{Strategy adaptation}\label{strategy-adaptation}}

In the adaptive ABCW model used here, each agent updates its strategy
using only the outcome of the immediately preceding game. A winner
retains its strategy, whereas a loser reverses it. Thus,

\({s_i(t+1)= \begin{cases} s_i(t),&u_i(t)=+1,\cr -s_i(t),&u_i(t)=-1 \end{cases}}\)

or, equivalently,

\[{s_i(t+1)=u_i(t)s_i(t).}\]

\hypertarget{action-update}{%
\subsubsection{Action update}\label{action-update}}

After payoffs have been calculated from the current actions and
reference signals from the current field, strategies are updated first.
A losing agent uses this \textbf{updated strategy}, rather than its
previous strategy, to determine its next action. The action update is
therefore

\({a_i(t+1)= \begin{cases} a_i(t),&u_i(t)=+1,\cr a_i(t),&u_i(t)=-1,\ \sigma_i(t)=0,\cr \operatorname{sgn}\!\left(s_i(t+1)\sigma_i(t)\right),&u_i(t)=-1,\ \sigma_i(t)\neq0 \end{cases}.}\)

A winner retains its current action. A loser also retains its current
action if no effective reference signal is available; otherwise, it
selects the next action according to its updated strategy and the
reference signal.

\hypertarget{field-update}{%
\subsubsection{Field update}\label{field-update}}

Finally, the payoff is fed back into the field. Let the learning rate be
\({\eta>0}\). All experiments analyzed in this paper use

\[{\eta=1.}\]

The field update is

\({w_{ij}(t+1)= \begin{cases} \max\!\left(0,w_{ij}(t)+\eta u_i(t)\right),&E_{ij}=1,\cr 0,&E_{ij}=0 \end{cases}.}\)

Thus, when agent \({i}\) wins, the weights of the existing outgoing
edges from that agent increase; when it loses, they decrease, subject to
the lower bound of zero.

The update order for one game is fixed as

\({\begin{array}{c} (a(t),s(t),W(t))\cr \downarrow\cr u(t),\sigma(t)\cr \downarrow\cr s(t+1)\cr \downarrow\cr a(t+1)\cr \downarrow\cr W(t+1) \end{array}.}\)

This ordering is part of the model definition because losing agents use
their updated, rather than previous, strategies when selecting their
next actions.

Every branch of the payoff, reference, strategy, action, and field
updates is resolved by these rules. In particular, tied reference agents
are averaged rather than selected arbitrarily, a zero reference signal
invokes the action-retention rule above, and no random tie-breaking or
stochastic choice is used. Consequently, for fixed
\({(a(t),s(t),W(t))}\) and \({E}\), the complete successor state is
unique. The ABCW system analyzed here is therefore deterministic, and
repeated runs from the same initial state generate the same finite
trajectory and transition records.

Figure 2 summarizes the complete update cycle used to generate the
dataset.

\begin{figure}
\centering
\includegraphics[trim=0 0 0 24bp,clip]{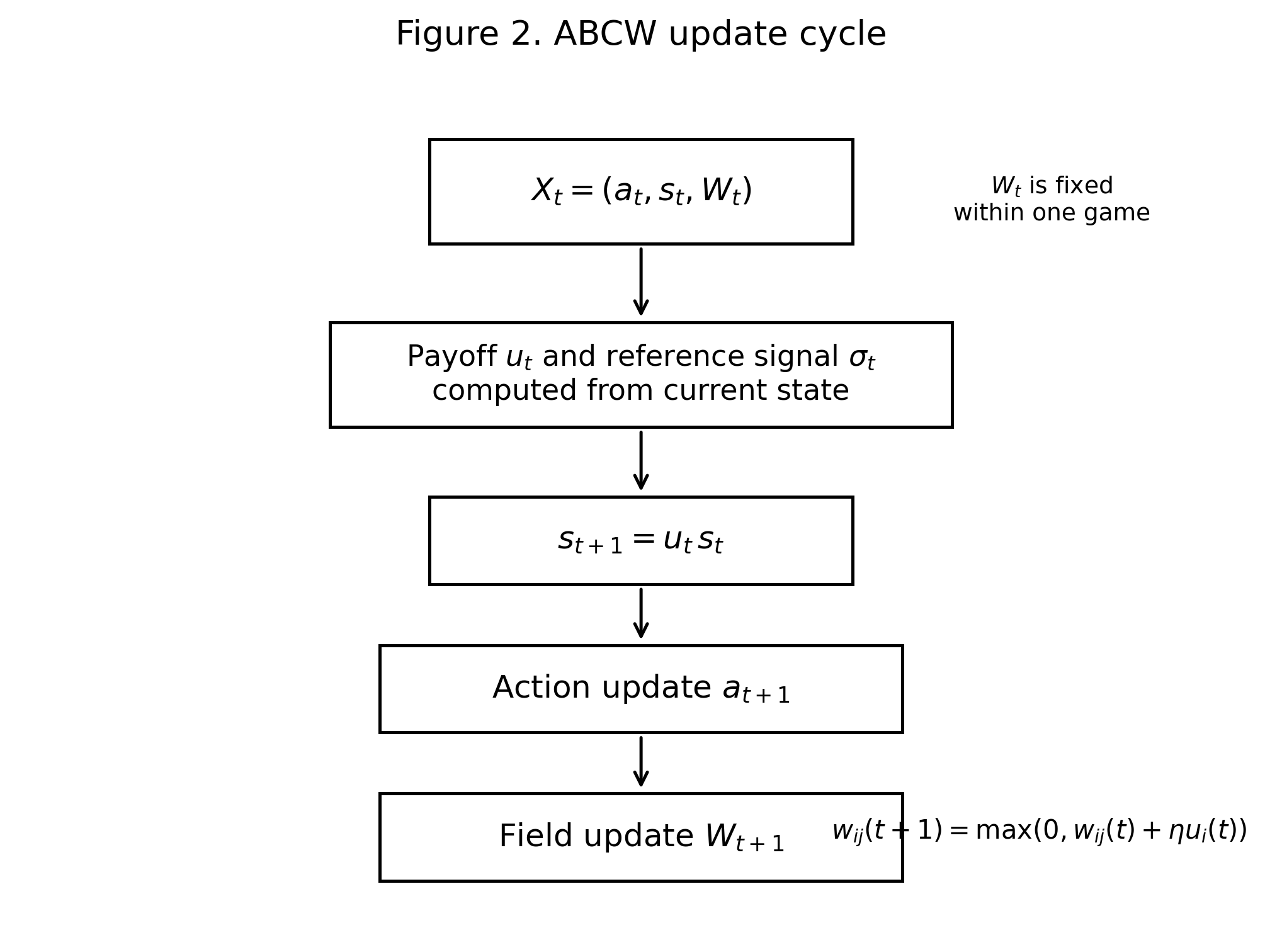}
\caption{ABCW update cycle. Update order for the adaptive
five-agent ABCW model used to construct the dataset. Payoffs and
reference signals are computed from the current state; strategies are
updated, then actions, then the influence field.}
\end{figure}

\hypertarget{field-representation}{%
\subsection{Field Representation}\label{field-representation}}

For the five-agent system, define the fully symmetric reference field as

\({W^\ast= \begin{array}{ccccc} 0&1&1&1&1\cr 1&0&1&1&1\cr 1&1&0&1&1\cr 1&1&1&0&1\cr 1&1&1&1&0 \end{array}.}\)

The reference field \({W^\ast}\) is not assumed to be an equilibrium or
a naturally privileged social state. It serves only as a reference point
for describing field deformation.

The field deviation is defined as

\[{\Delta W(t)=W(t)-W^\ast,}\]

or componentwise as

\[{\Delta w_{ij}(t)=w_{ij}(t)-w^\ast_{ij}.}\]

Thus, \({\Delta W(t)=0}\) means that the current field coincides with
the reference field. A nonzero deviation \({\Delta W\neq0}\) does not
necessarily imply matrix asymmetry \({W\neq W^{\mathsf T}}\). Throughout
this paper, \emph{field deviation} refers to the difference from the
reference field itself. This \({\Delta W}\) is the field variable used
in the subsequent prediction problem.

\hypertarget{dataset-construction}{%
\subsection{Dataset Construction}\label{dataset-construction}}

The analysis uses four initial fields for the five-agent ABCW model:

\begin{itemize}
\tightlist
\item
  Baseline
\item
  Hub
\item
  Local
\item
  Hub+Local
\end{itemize}

For Baseline, \({W(0)=W^\ast}\). In each of the other three initial
fields, four directed edges whose reference weight is 1 are strengthened
to 4, while the update rules remain unchanged. Only the placement of
these strong edges differs. Numbering the agents \({1,\ldots,5}\) and
denoting the directed edge \({i\to j}\) by \({(i,j)}\), the strengthened
edges are:

\begin{itemize}
\tightlist
\item
  Hub: \({(1,2),(1,3),(1,4),(1,5)}\)
\item
  Local: \({(1,2),(2,3),(3,4),(4,5)}\)
\item
  Hub+Local: \({(1,3),(1,5),(2,3),(4,5)}\)
\end{itemize}

In Hub, the four outgoing edges from agent 1 to the other agents are
strengthened. Local strengthens a chain of four edges. Hub+Local
combines two outgoing edges from agent 1 with two local competing edges
that share their destination agents.

Each of these three non-baseline fields therefore has four nonzero
initial deviation components and the same Frobenius norm,

\[{\lVert\Delta W(0)\rVert_F=\sqrt{4\times3^2}=6.}\]

These configurations are not intended as an exhaustive classification of
network structures. They serve as initial-condition families from which
a finite reachable state set is constructed.

For each initial field, we exhaust all combinations of the five-agent
initial action vector \({a(0)}\) and initial strategy vector \({s(0)}\).
Because every component is binary, this gives

\[{2^5\times2^5=1024}\]

initial conditions per field and

\[{4\times1024=4096}\]

initial conditions in total.

Each case is evolved until the complete state recurs, and every
transition encountered before recurrence is collected. This procedure
yields 56,536 transitions containing 2,562 distinct current fields
\({\Delta W}\) and 11,202 distinct \({(a,\Delta W)}\) states.

Table 2 summarizes the dataset construction and the resulting counts.

\hypertarget{table-2.-dataset-construction-and-summary}{%
\subsection{Table 2. Dataset construction and
summary}\label{table-2.-dataset-construction-and-summary}}

\begin{longtable}[]{@{}lr@{}}
\toprule\noalign{}
Quantity & Value \\
\midrule\noalign{}
\endhead
\bottomrule\noalign{}
\endlastfoot
Agents & 5 \\
Learning rate \(\eta\) & 1 \\
Initial field configurations & 4 \\
Action configurations per field & 32 \\
Strategy configurations per field & 32 \\
Initial conditions & 4096 \\
Observed transitions & 56536 \\
Distinct current fields & 2562 \\
Distinct \((a,\Delta W)\) states & 11202 \\
\end{longtable}

The 56,536 transitions do not represent the complete theoretical state
space of ABCW. All subsequent results should be interpreted as
statements about the stored finite transition set reached from the four
initial fields and 4,096 initial conditions defined above. We make no
claim of generalization to unreachable states, other numbers of agents,
or other initial fields. Nor do we interpret the recurrence of every
trajectory among the 4,096 cases as a proof that the weights are bounded
in general.

\begin{center}\rule{0.5\linewidth}{0.5pt}\end{center}

\hypertarget{observation-and-predictive-objective}{%
\section{Observation and Predictive
Objective}\label{observation-and-predictive-objective}}

This section formulates, for the finite ABCW system defined above, the
problem of determining which distinctions in the present must be
retained to predict a specified future exactly. The central point is to
distinguish closure of the observed dynamics from sufficiency for a
particular prediction target.

\hypertarget{observation-as-information-loss}{%
\subsection{Observation as Information
Loss}\label{observation-as-information-loss}}

Let \({\Omega}\) denote the complete state space, and represent an
observation as a mapping

\[{P:\Omega\longrightarrow Y.}\]

If two distinct complete states \({X,X'\in\Omega}\) satisfy

\[{P(X)=P(X'),}\]

then the observation \({P}\) discards the distinction between them and
treats them as the same observed state. An observation is therefore not
merely a choice of which variables to inspect. It is a mapping that
determines which distinctions among complete states are retained and
which are identified.

The complete ABCW state is

\[{X=(a,s,W).}\]

Because the reference field \({W^\ast}\) is fixed, \({W}\) and
\({\Delta W}\) are in one-to-one correspondence. For example,

\[{P_W(X)=\Delta W}\]

discards distinctions in action \({a}\) and strategy \({s}\), whereas

\[{P_{aW}(X)=(a,\Delta W)}\]

discards only distinctions in strategy \({s}\).

\hypertarget{self-closure-and-target-prediction}{%
\subsection{Self-Closure and Target
Prediction}\label{self-closure-and-target-prediction}}

Write the deterministic ABCW update as

\[{X_{t+1}=\Phi(X_t).}\]

An observation \({P}\) is called \textbf{self-closed} if

\[{P(X)=P(X')}\]

always implies

\[{P(\Phi(X))=P(\Phi(X')).}\]

Under this condition, the next observed state is uniquely determined by
the current observed state. Self-closure is a deterministic
quotient-consistency condition that makes the coarse-grained transition
well defined. It is related to aggregation conditions in lumpability and
bisimulation, but it is not identical to causal-state construction from
conditional distributions over future sequences
(\protect\hyperlink{ref-shalizi2001}{Shalizi and Crutchfield 2001};
\protect\hyperlink{ref-givan2003}{Givan, Dean, and Greig 2003}).

Self-closure is distinct from sufficiency for a specified prediction
target. If the quantity of interest is

\[{Z=g(X_{t+1}),}\]

the required condition is that

\[{P(X_t)\longrightarrow Z}\]

be single-valued. It is not necessary for

\[{P(X_t)\longrightarrow P(X_{t+1})}\]

to be single-valued.

This distinction appears directly in the finite ABCW dataset. For the
field-only observation

\[{P_W(X)=\Delta W,}\]

the next field is unique for 1,390 of the 2,562 observed current fields:

\[{\frac{1390}{2562}\simeq0.54254.}\]

Thus, the current field alone determines the complete next field for
approximately 54.3\% of the observed current fields, but not in general.

When the current action is retained through

\[{P_{aW}(X)=(a,\Delta W),}\]

the relation

\[{(a_t,\Delta W_t)\longrightarrow\Delta W_{t+1}}\]

is single-valued for all 11,202 observed \({(a,\Delta W)}\) states. By
contrast, the next state of this observation itself,

\[{(a_t,\Delta W_t)\longrightarrow(a_{t+1},\Delta W_{t+1}),}\]

is single-valued for only 8,390 of the 11,202 observed states:

\[{\frac{8390}{11202}\simeq0.748973.}\]

This is approximately 74.9\%. An observation that is not self-closed can
therefore still predict a specified component of the future perfectly.
We keep these two requirements separate and fix the prediction target
explicitly below.

\hypertarget{one-step-exact-field-prediction}{%
\subsection{One-Step Exact-Field
Prediction}\label{one-step-exact-field-prediction}}

The prediction target in this paper is

\[{Z=\Delta W_{t+1}.}\]

That is, the observation must predict the \textbf{complete
one-step-ahead field deviation} without error.

Retaining \({(a_t,\Delta W_t)}\) makes this target single-valued
throughout the finite dataset, but it also retains all 2,562 distinct
field deviations. We therefore keep the current action \({a_t}\)
unchanged while further coarsening only the field component.

Let

\[{B:\Delta W\longrightarrow\mathcal C}\]

be a field-only map, and define the corresponding observation as

\[{P_B(X_t)=\left(a_t,B(\Delta W_t)\right).}\]

The map \({B}\) satisfies the predictive requirement if, for any two
current states \({X_t}\) and \({X'_t}\) in the analyzed transitions,

\[{a_t=a'_t}\]

and

\[{B(\Delta W_t)=B(\Delta W'_t)}\]

imply

\[{\Delta W_{t+1}=\Delta W'_{t+1}.}\]

Equivalently, we require

\[{\left(a_t,B(\Delta W_t)\right)\longrightarrow\Delta W_{t+1}}\]

to be single-valued on the stored finite transition set.

The action \({a}\) is not itself compressed, and the prediction target
is not coarse-grained: the required output remains the complete next
field deviation. The problem considered here is therefore an
\textbf{action-conditioned, one-step, exact-output, field-only
predictive-partition problem}.

We call the required single-valuedness \textbf{predictive sufficiency}
for the specified target. In the present formulation, \({B}\) is
admissible when \({(a_t,B(\Delta W_t))}\) is predictively sufficient for
\({\Delta W_{t+1}}\) on the analyzed dataset.

\hypertarget{minimum-field-partition}{%
\subsection{Minimum Field Partition}\label{minimum-field-partition}}

Field deviations mapped to the same value by \({B}\) are treated as
members of the same class, and the set of classes is
\({\operatorname{Im}B}\). The complete field observation has

\[{|\operatorname{Im}B|=2562.}\]

If exact predictive accuracy can be retained with

\[{|\operatorname{Im}B|<2562,}\]

then distinct field deviations can be identified for the specified
prediction task, at least on the finite reachable set considered here.

We define the optimization problem as

\[{\min_B|\operatorname{Im}B|}\]

subject to

\({\left(a_t,B(\Delta W_t)\right)\longrightarrow\Delta W_{t+1} \quad\text{being single-valued on all observed transitions}.}\)

This formulation is related in motivation to the Information Bottleneck,
which compresses a representation while retaining relevant information.
It is not the standard Information Bottleneck objective. Rather than
optimizing a trade-off between compression and relevant information
through mutual information, we impose zero prediction error as a hard
constraint and minimize the finite-set class count
\({|\operatorname{Im}B|}\) (\protect\hyperlink{ref-tishby2000}{Tishby,
Pereira, and Bialek 2000}).

The problem also differs from finding the coarsest self-closed
partition. Here, the retained action serves as a conditioning variable,
and the quantity preserved is the external target \({\Delta W_{t+1}}\).
The minimization is therefore defined by target prediction rather than
self-closure of the observation.

When useful, we report the field-class compression rate

\[{C_{\Delta W}=1-\frac{|\operatorname{Im}B|}{2562}.}\]

This is a class-count-based compression rate. It is not a reduction rate
in Shannon information or in the number of bits required for
representation.

This formulation separates the minimization problem from the prior
selection of features that appear intuitively natural. Section 5 first
evaluates baselines constructed from natural network features. Section 6
then derives the required field distinctions directly from the
predictive constraint and determines the minimum partition.

\begin{center}\rule{0.5\linewidth}{0.5pt}\end{center}

\hypertarget{natural-feature-baselines}{%
\section{Natural-Feature Baselines}\label{natural-feature-baselines}}

Before deriving the minimum predictive partition directly from the
dynamics, we examine how coarsely the current field \({\Delta W}\) can
be represented by human-selected network features, including norms,
local structure, hub structure, and out-strength. The purpose is to
assess whether interpretability, genuine field compression, and
predictive sufficiency for the complete one-step-ahead field can be
achieved simultaneously.

\hypertarget{natural-feature-compression-and-its-limit}{%
\subsection{Natural-feature compression and its
limit}\label{natural-feature-compression-and-its-limit}}

The out-strength vector, for example, achieves very high predictive
performance. Its state-level determinism score is evaluated over the
observed states defined by \({(a,\mathrm{out\mbox{-}strength})}\). The
next field is unique for 11,084 of the 11,142 distinct observed states:

\[{\frac{11084}{11142}=0.994794\ldots.}\]

The resulting state-level determinism score is 99.4794\%. The remaining
58 observed states are conflicts: more than one next field occurs from
the same observed state.

The denominator here is not the total number of 56,536 transitions but
the 11,142 distinct observed states. This state-level score should
therefore be distinguished from evaluations weighted over the complete
transition set.

We then searched for simple feature combinations that would resolve
these 58 conflicts while genuinely compressing \({\Delta W}\). Within
the candidate family examined, no feature representation simultaneously
achieved nontrivial field compression and exact predictive determinism.

This result contrasts with the response-based minimum predictive
partition constructed exactly in Section 6. That partition compresses
2,562 fields into 692 classes while preserving exact one-step field
prediction across all 56,536 observed transitions.

The contrast reflects two different approaches to compression. The first
summarizes the field using features selected in advance:

\[{\Delta W\longrightarrow B(\Delta W).}\]

The second derives the necessary distinctions from equivalence of future
responses:

\[{\Delta W\longrightarrow\left(a\longmapsto\Delta W'\right).}\]

The first approach favors interpretability but may discard information
required for prediction. The second preserves distinctions required by
the predictive task, but the resulting classes need not have a simple
interpretation as network statistics. The 692-class partition derived in
Section 6 belongs to the second category.

\hypertarget{exhaustive-combinations-of-the-candidate-feature-families}{%
\subsection{Exhaustive combinations of the candidate feature
families}\label{exhaustive-combinations-of-the-candidate-feature-families}}

The exhaustive search uses the following nine feature maps. For a field
matrix \({D=\Delta W}\) with singular values
\({\sigma_1(D)\ge\cdots\ge\sigma_5(D)\ge0}\), define its ordered out-
and in-strength vectors by

\[{r(D)=D\mathbf 1,\qquad c(D)=D^{\mathsf T}\mathbf 1.}\]

The nine maps are:

\begin{enumerate}
\def\labelenumi{\arabic{enumi}.}
\tightlist
\item
  \textbf{Squared Frobenius norm:}
  \({\lVert D\rVert_F^2=\sum_{i,j}D_{ij}^2}\).
\item
  \textbf{Matrix rank:} \({\operatorname{rank}(D)}\), evaluated
  numerically on the integer-valued matrix.
\item
  \textbf{Stable rank:}
  \({\operatorname{srank}(D)=\lVert D\rVert_F^2/\sigma_1(D)^2}\), with
  value zero when \({\sigma_1(D)=0}\).
\item
  \textbf{Singular-value vector:}
  \({(\sigma_1(D),\ldots,\sigma_5(D))}\).
\item
  \textbf{Out-strength distribution:} the entries of \({r(D)}\) sorted
  in nondecreasing order.
\item
  \textbf{In-strength distribution:} the entries of \({c(D)}\) sorted in
  nondecreasing order.
\item
  \textbf{Joint in/out-strength distribution:} the five pairs
  \({(c_i(D),r_i(D))}\) sorted lexicographically.
\item
  \textbf{Out-strength vector:} the ordered vector \({r(D)}\), retaining
  agent identities.
\item
  \textbf{In-strength vector:} the ordered vector \({c(D)}\), retaining
  agent identities.
\end{enumerate}

Thus, ``distribution'' denotes a permutation-invariant sorted
representation, whereas ``vector'' preserves the ordering of the five
agents. Numerically computed singular values and stable rank are
canonicalized to ten significant digits in the exhaustive comparison;
the remaining feature values are integer-valued. A feature combination
is the tuple of its selected map values, with no learned weighting or
post-processing.

These nine feature maps have

\[{2^9-1=511}\]

nonempty subsets. We exhaustively evaluated all 511 feature
combinations. None simultaneously merged genuinely distinct values of
\({\Delta W}\) and determined the complete one-step-ahead field without
ambiguity.

The out-strength vector yields high determinism but leaves 58
conflicting observed states. A feature representation denoted by
\({B^\ast}\) achieves exact prediction, but it uniquely distinguishes
all 2,562 fields and therefore provides no genuine field compression.

Table 3 compares representative natural-feature baselines with the exact
minimum partition. Here \({D_{state}}\) is the state-level determinism
score, whereas \({A_{freq}}\) is the corresponding
transition-frequency-weighted score.

\hypertarget{table-3.-natural-feature-baselines-and-exact-minimum-partition}{%
\subsection{Table 3. Natural-feature baselines and exact minimum
partition}\label{table-3.-natural-feature-baselines-and-exact-minimum-partition}}

\begin{longtable}[]{@{}
  >{\raggedright\arraybackslash}p{(\columnwidth - 10\tabcolsep) * \real{0.1667}}
  >{\raggedleft\arraybackslash}p{(\columnwidth - 10\tabcolsep) * \real{0.1667}}
  >{\raggedleft\arraybackslash}p{(\columnwidth - 10\tabcolsep) * \real{0.1667}}
  >{\raggedleft\arraybackslash}p{(\columnwidth - 10\tabcolsep) * \real{0.1667}}
  >{\raggedleft\arraybackslash}p{(\columnwidth - 10\tabcolsep) * \real{0.1667}}
  >{\raggedright\arraybackslash}p{(\columnwidth - 10\tabcolsep) * \real{0.1667}}@{}}
\toprule\noalign{}
\begin{minipage}[b]{\linewidth}\raggedright
Observation
\end{minipage} & \begin{minipage}[b]{\linewidth}\raggedleft
Field classes
\end{minipage} & \begin{minipage}[b]{\linewidth}\raggedleft
Observed \((a,B)\) states
\end{minipage} & \begin{minipage}[b]{\linewidth}\raggedleft
\(D_{state}\)
\end{minipage} & \begin{minipage}[b]{\linewidth}\raggedleft
\(A_{freq}\)
\end{minipage} & \begin{minipage}[b]{\linewidth}\raggedright
True field compression
\end{minipage} \\
\midrule\noalign{}
\endhead
\bottomrule\noalign{}
\endlastfoot
Full field \(\Delta W\) & 2562 & 11202 & 1 & 1 & No \\
Out-strength vector & 2516 & 11142 & 0.994794 & 0.994658 & Yes \\
Best true-compression candidate by \(D_{state}\) (in-strength
distribution + out-strength vector) & 2561 & 11200 & 0.999821 & 0.999717
& Yes \\
\(B^\ast\) (Frobenius norm + out-strength vector)
& 2562 & 11202 & 1 & 1 & No \\
Exact minimum partition \(B_{692}\) & 692 & 8084 & 1 & 1 & Yes \\
\end{longtable}

Figure 3 places all 511 candidates in the same performance--compression
plane and highlights the exact minimum partition.

\begin{figure}
\centering
\includegraphics[trim=0 0 0 24bp,clip]{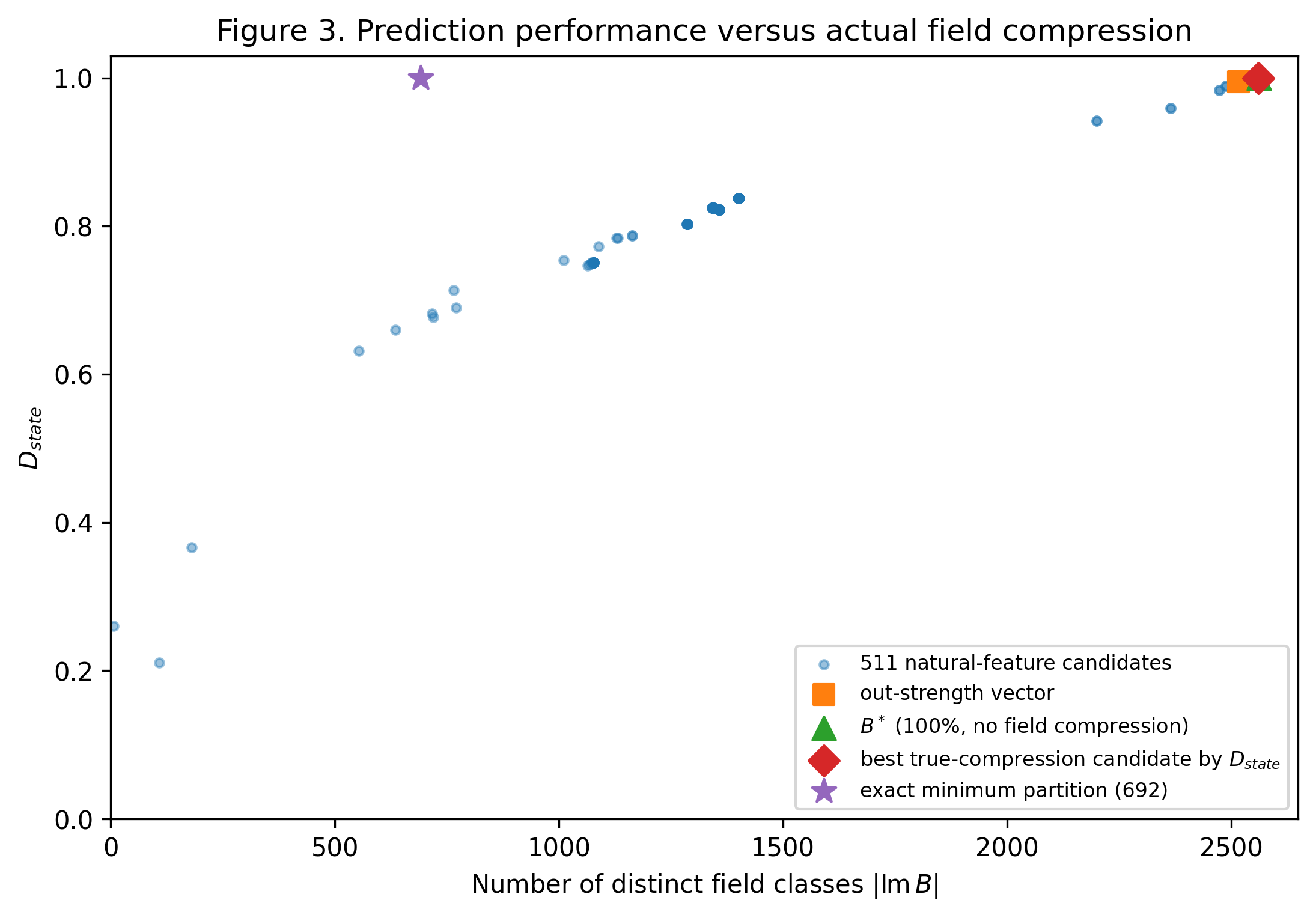}
\caption{Prediction performance versus actual field
compression. All 511 nonempty combinations from the nine natural feature
families are plotted using the number of distinct field values and the
state-level determinism score. The out-strength vector, the
100\%-predictive but noncompressive B* representation, the best
genuinely compressive candidate by D\_state, and the exact 692-class
partition are highlighted.}
\end{figure}

The failure of the natural-feature search to achieve both exact
prediction and genuine compression does not imply that nearly all 2,562
current fields must remain distinct. It shows only that the tested
feature family does not express the required partition. In the next
section, we remove the restriction of selecting features in advance and
derive the necessary distinctions directly from conflicts in future
responses.

\begin{center}\rule{0.5\linewidth}{0.5pt}\end{center}

\hypertarget{exact-minimum-predictive-partition}{%
\section{Exact Minimum Predictive
Partition}\label{exact-minimum-predictive-partition}}

\hypertarget{incompatibility-graph}{%
\subsection{Incompatibility graph}\label{incompatibility-graph}}

We now remove the restriction of preselected features and derive
directly from the dynamics which distinctions must be retained to avoid
confusing future fields. The analyzed finite reachable set is the same
as in Sections 3--5 and contains:

\begin{itemize}
\tightlist
\item
  56,536 transitions;
\item
  2,562 distinct current fields; and
\item
  11,202 distinct \({(a,\Delta W)}\) states.
\end{itemize}

Consider a field-only map

\[{B:\Delta W\longrightarrow\mathcal C}\]

and the observation

\[{P_B(X_t)=(a_t,B(\Delta W_t)).}\]

The requirement is that

\[{(a_t,B(\Delta W_t))\longrightarrow\Delta W_{t+1}}\]

be single-valued over all 56,536 observed transitions. In other words,
the same action and the same field class must never produce different
next fields.

Consider two current fields \({\Delta W_i}\) and \({\Delta W_j}\). If
both are observed under a common action \({a}\) and produce different
next fields \({\Delta W'_i\neq\Delta W'_j}\), they cannot be assigned to
the same field class. If

\[{B(\Delta W_i)=B(\Delta W_j),}\]

then the identical observation

\[{(a,B(\Delta W_i))=(a,B(\Delta W_j))}\]

would lead to two different next fields, violating single-valued exact
prediction.

We therefore construct an incompatibility graph \({G}\) with the 2,562
distinct current fields as vertices. An edge joins two fields if they
produce different next fields under at least one common observed action.
A field-only partition preserving exact prediction is then equivalent to
a proper coloring of \({G}\), because adjacent vertices must receive
different class labels.

It follows that the minimum number of field classes is

\[{|\operatorname{Im}B_{\min}|=\chi(G),}\]

where \({\chi(G)}\) is the chromatic number of the incompatibility
graph.

This formulation does not introduce a new general state-minimization
problem. It is connected to compatibility-based minimization of
partially specified input-output behavior, particularly a depth-1
partial Mealy-machine formulation of incompletely specified finite-state
machine minimization. The objective here is to obtain the concrete exact
solution produced when this established problem form is applied to the
finite reachable dynamics of ABCW.

\hypertarget{constructive-upper-bound}{%
\subsection{Constructive upper bound}\label{constructive-upper-bound}}

We first color the incompatibility graph using DSATUR
(\protect\hyperlink{ref-brelaz1979}{Brélaz 1979}). This produces a valid
proper coloring with 692 colors and hence establishes

\[{\chi(G)\le692.}\]

The number returned by a coloring algorithm is not, by itself,
sufficient to verify the original prediction requirement. We therefore
interpret the resulting coloring as the field-only map

\[{B_{692}(\Delta W)}\]

and map it back to all 56,536 original transitions. We then directly
test whether

\[{(a,B_{692}(\Delta W))\longrightarrow\Delta W_{t+1}}\]

is single-valued. The result is:

\begin{itemize}
\tightlist
\item
  field classes: 692;
\item
  distinct \({(a,B_{692})}\) states: 8,084;
\item
  nondeterministic observed states: 0; and
\item
  one-step complete-field prediction: 100\%.
\end{itemize}

Thus, a field-only partition exists that compresses

\[{2562\longrightarrow692}\]

while preserving exact one-step prediction of the complete field. The
class-count-based compression rate is

\[{C_{\Delta W}=1-\frac{692}{2562}\simeq0.7299,}\]

or approximately 72.99\%. The 11,202 complete \({(a,\Delta W)}\)
observations are reduced to 8,084 distinct \({(a,B_{692}(\Delta W))}\)
observations.

Therefore, 692 is an upper bound:

\[{\boxed{\chi(G)\le692}.}\]

At this stage, however, a proper coloring with 691 or fewer colors has
not been ruled out. Establishing exact minimality requires an
independent lower bound of at least 692.

\hypertarget{anchor-lower-bound}{%
\subsection{Anchor lower bound}\label{anchor-lower-bound}}

To construct a lower bound, we fix one action condition. We choose the
anchor action

\[{a=(-,-,-,-,-).}\]

Under this action, 1,239 current fields are observed, and they produce
625 distinct next fields. The lower-bound argument does not assume that
these anchor groups cover all 2,562 fields. It uses only the induced
subgraph on the 1,239 vertices actually observed under the anchor
action.

Two current fields producing different next fields under the same anchor
action cannot belong to the same field class. Partition the 1,239 fields
into 625 anchor future groups

\[{V_1,\ldots,V_{625}}\]

according to their next field under the anchor action. Every pair of
vertices belonging to different groups is adjacent in \({G}\).

\hypertarget{lemma-1-additivity-across-anchor-groups}{%
\subsubsection{Lemma 1 --- Additivity across anchor
groups}\label{lemma-1-additivity-across-anchor-groups}}

For two distinct anchor future groups \({V_i}\) and \({V_j}\), with
\({i\neq j}\), let

\[{u\in V_i,\qquad v\in V_j.}\]

By construction, \({u}\) and \({v}\) produce different next fields under
the anchor action. Therefore,

\[{\lbrace u,v\rbrace\in E(G),}\]

and every two distinct anchor future groups form a complete join.
Consequently, different anchor groups cannot share a color in any proper
coloring.

Let

\[{G_i=G[V_i]}\]

be the subgraph induced by group \({V_i}\). The chromatic number of the
anchor-observed induced subgraph is then additive:

\({\chi\!\left(G\!\left[\bigcup_{i=1}^{625}V_i\right]\right) =\sum_{i=1}^{625}\chi(G_i).}\)

Because every group requires at least one color,

\[{\chi(G)\ge625.}\]

Combining this with the constructive upper bound gives

\[{\boxed{625\le\chi(G)\le692}.}\]

The value 625 captures only the distinctions exposed by the anchor
action. Two fields that lead to the same anchor future may still lead to
different futures under another action condition. Such fields must be
separated when all observed actions are considered simultaneously.

We therefore inspect the incompatibilities remaining within each of the
625 anchor future groups. The result is:

\begin{itemize}
\tightlist
\item
  575 groups with no internal conflicts; and
\item
  50 groups with internal conflicts.
\end{itemize}

For the 575 conflict-free groups, \({\chi(G_i)=1}\). The chromatic
numbers of the remaining 50 induced subgraphs must be determined
exactly.

\begin{figure}
\centering
\includegraphics[trim=0 0 0 24bp,clip]{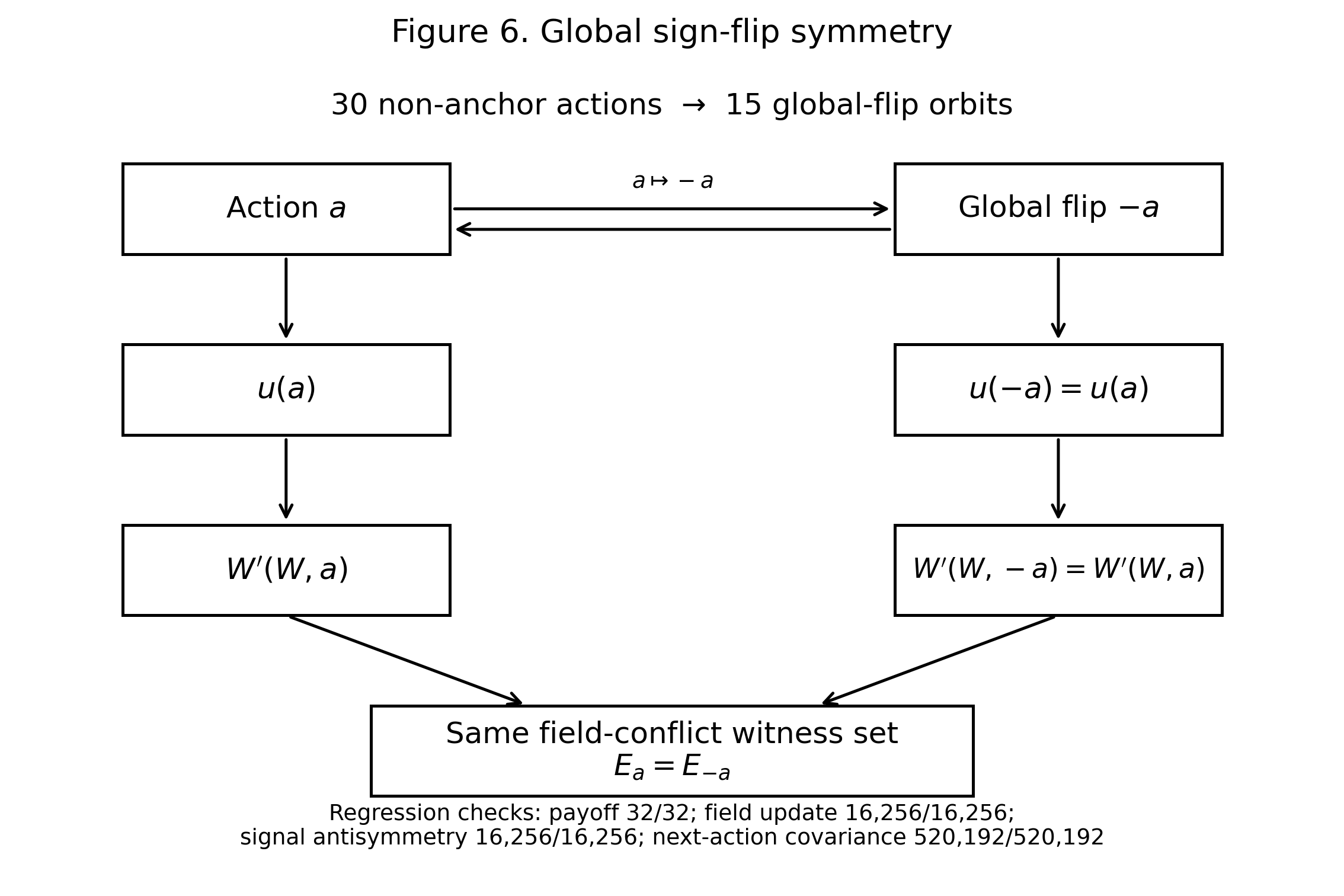}
\caption{Global sign-flip symmetry. Under the global action
flip a -\textgreater{} -a, minority-game payoffs are unchanged and the
one-step field update is invariant, so a and -a witness exactly the same
field-conflict edges. Hence 30 non-anchor actions collapse to 15
global-flip orbits.}
\end{figure}

Under this symmetry, actions \({a}\) and \({-a}\) witness the same set
of field-conflict edges. The 30 non-anchor actions can therefore be
grouped into 15 global-flip orbits for the purpose of describing
internal conflicts. This symmetry does not itself provide the lower
bound of 692. It identifies duplicated conflict information and supplies
both a redundancy reduction and a consistency check for interpreting the
internal analysis of the 50 difficult groups.

\hypertarget{exact-refinement-of-the-50-difficult-anchor-groups}{%
\subsection{Exact refinement of the 50 difficult anchor
groups}\label{exact-refinement-of-the-50-difficult-anchor-groups}}

For each of the 50 anchor groups containing internal conflicts, we
compute the chromatic number of the induced incompatibility graph
\({G_i}\) exactly.

For each group, a greedy DSATUR coloring provides an upper bound
(\protect\hyperlink{ref-brelaz1979}{Brélaz 1979}), and a
Bron--Kerbosch-type maximum-clique search provides a lower bound
(\protect\hyperlink{ref-bron1973}{Bron and Kerbosch 1973}). If the
bounds coincide, their common value determines \({\chi(G_i)}\)
immediately. If they do not, we test \({k}\)-colorability by
backtracking in DSATUR order for successive values of \({k}\) from the
lower bound to the upper bound. The smallest feasible value is the exact
chromatic number.

The values reported below are therefore not merely the color counts
returned by a greedy heuristic. For every induced subgraph, they are
certified by matching lower bounds and colorability results.

Across all 625 anchor groups, the exact chromatic-number distribution
is:

\begin{itemize}
\tightlist
\item
  \({\chi=1}\): 575 groups;
\item
  \({\chi=2}\): 39 groups;
\item
  \({\chi=3}\): 8 groups;
\item
  \({\chi=4}\): 2 groups; and
\item
  \({\chi=7}\): 1 group.
\end{itemize}

The group requiring the largest number of colors, \({\chi=7}\), has 44
vertices and 224 internal conflict edges. Its maximum-clique lower bound
and DSATUR upper bound are both 7, directly establishing

\[{\chi(G_i)=7.}\]

The vertex count, internal edge count, exact chromatic number, lower
bound, and upper bound for every group are stored in the reproducibility
data.

By Lemma 1, colors cannot be shared across distinct anchor future
groups. The total number of colors required by the anchor-observed
induced subgraph is therefore the sum of the group chromatic numbers.
Relative to assigning one color to each of the 625 groups, the
additional number of colors is

\[{39(2-1)+8(3-1)+2(4-1)+1(7-1)}\]

\[{=39+16+6+6}\]

\[{=67.}\]

Hence,

\[{625+67=692.}\]

The same result is obtained by the direct cross-check

\[{575\cdot1+39\cdot2+8\cdot3+2\cdot4+1\cdot7=692.}\]

Thus, the induced subgraph on the 1,239 fields observed under the anchor
action already requires 692 colors. The complete incompatibility graph
must therefore satisfy

\[{\boxed{\chi(G)\ge692}.}\]

The value 692 is consequently not an accidental output of DSATUR. The
anchor action alone forces 625 distinctions. The remaining 67 arise
because some fields that are indistinguishable under the anchor
condition must be separated under other action conditions. Moreover,
these additional distinctions are localized in only 50 of the 625
groups.

The decomposition

\[{692=625+67}\]

therefore reflects the complete-join structure across anchor groups and
the action-conditioned conflicts remaining within them. Figure 5
summarizes the matching upper and lower bounds and this decomposition.

Table 4 records the exact contribution of each anchor-group chromatic
number to the lower bound.

\hypertarget{table-4.-exact-lower-bound-decomposition}{%
\subsection{Table 4. Exact lower-bound
decomposition}\label{table-4.-exact-lower-bound-decomposition}}

\begin{longtable}[]{@{}
  >{\raggedright\arraybackslash}p{(\columnwidth - 6\tabcolsep) * \real{0.2500}}
  >{\raggedleft\arraybackslash}p{(\columnwidth - 6\tabcolsep) * \real{0.2500}}
  >{\raggedleft\arraybackslash}p{(\columnwidth - 6\tabcolsep) * \real{0.2500}}
  >{\raggedleft\arraybackslash}p{(\columnwidth - 6\tabcolsep) * \real{0.2500}}@{}}
\toprule\noalign{}
\begin{minipage}[b]{\linewidth}\raggedright
Required colors per anchor group
\end{minipage} & \begin{minipage}[b]{\linewidth}\raggedleft
Number of anchor groups
\end{minipage} & \begin{minipage}[b]{\linewidth}\raggedleft
Contribution to total classes
\end{minipage} & \begin{minipage}[b]{\linewidth}\raggedleft
Additional colors beyond one per group
\end{minipage} \\
\midrule\noalign{}
\endhead
\bottomrule\noalign{}
\endlastfoot
1 & 575 & 575 & 0 \\
2 & 39 & 78 & 39 \\
3 & 8 & 24 & 16 \\
4 & 2 & 8 & 6 \\
7 & 1 & 7 & 6 \\
Total & 625 & 692 & 67 \\
\end{longtable}

\begin{figure}
\centering
\includegraphics[trim=0 0 0 24bp,clip]{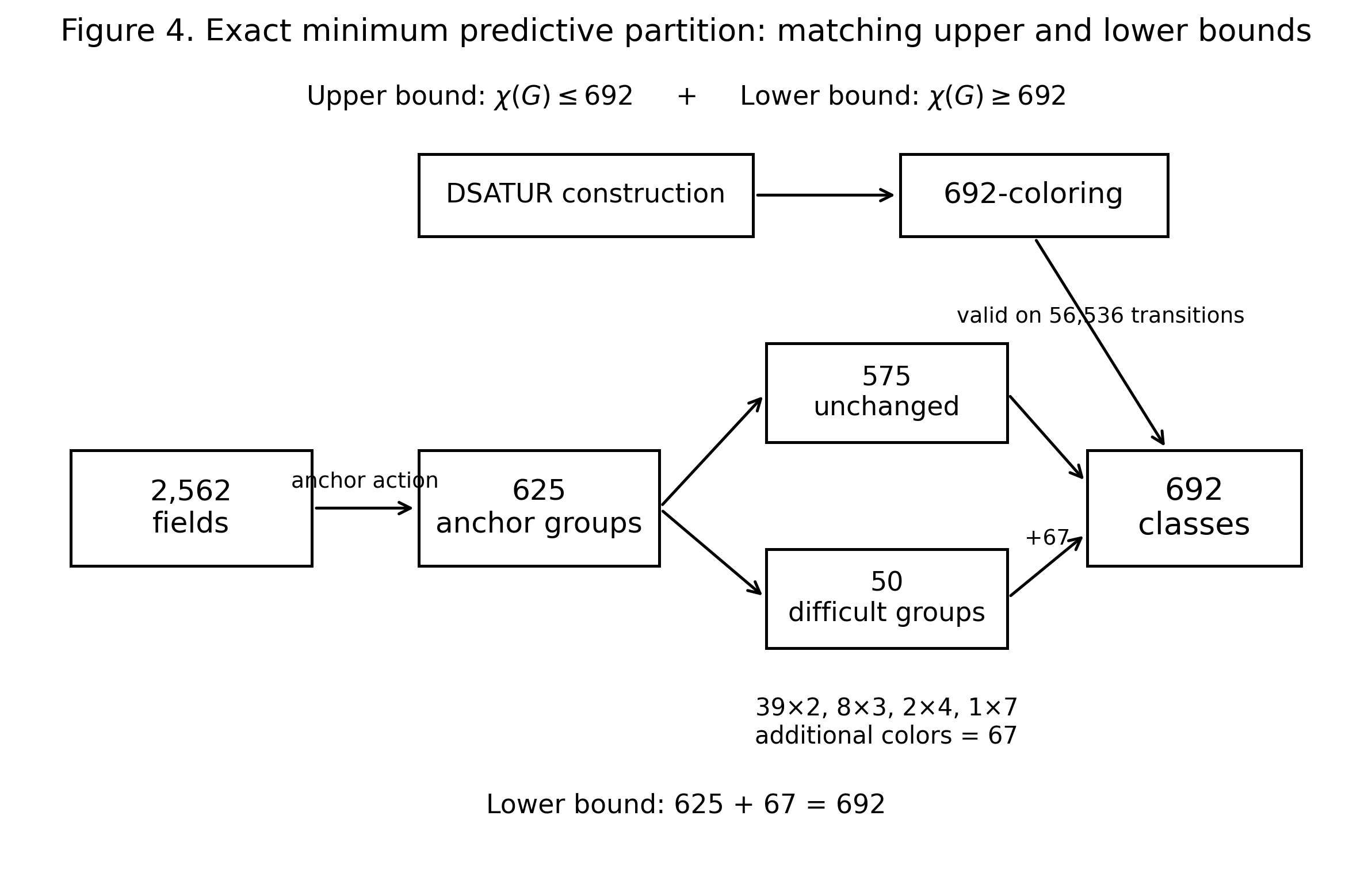}
\caption{Exact minimum predictive partition. The constructive
692-coloring provides the upper bound. The anchor action yields 625
groups; exact coloring of the 50 internally conflicting groups requires
67 additional colors, producing a matching lower bound of 692.}
\end{figure}

\hypertarget{main-result}{%
\subsection{Main result}\label{main-result}}

Section 6.2 constructs a valid 692-coloring of the incompatibility graph
on all 2,562 fields, establishing

\[{\chi(G)\le692.}\]

Direct evaluation of this coloring on all 56,536 original transitions
verifies that

\[{(a,B_{692}(\Delta W))\longrightarrow\Delta W_{t+1}}\]

is single-valued throughout the observed dataset.

Sections 6.3 and 6.4 decompose the 1,239 fields observed under the
anchor action into 625 future groups. The complete joins across
different groups and the exact chromatic numbers within the groups
establish

\[{\chi(G)\ge692.}\]

The bounds coincide:

\[{\boxed{\chi(G)=692}.}\]

\begin{quote}
\textbf{Main Result.}\\
For the 56,536 observed transitions and 2,562 distinct current fields of
the five-agent ABCW system, consider a field-only map

\({B:\Delta W\longrightarrow\mathcal C}\)

such that

\({(a_t,B(\Delta W_t))\longrightarrow\Delta W_{t+1}}\)

is single-valued over all observed transitions. The exact minimum number
of field classes is

\({\boxed{|\operatorname{Im}B_{\min}|=692}.}\)
\end{quote}

It is therefore unnecessary to retain all 2,562 distinct current fields
for the specified task. Exact prediction of the complete one-step-ahead
field permits the reduction

\[{2562\longrightarrow692,}\]

with the class-count-based compression rate

\[{1-\frac{692}{2562}\simeq72.99\%.}\]

Any partition into 691 or fewer classes fails even on the anchor-induced
subgraph used in the lower-bound argument. It must assign the same class
to at least one pair of fields that produce different next fields under
a common action condition, thereby violating exact prediction.

Thus, 692 is not merely the size of a high-performing candidate
partition. It is the minimum number of classes in a field-only partition
that preserves exact one-step, action-conditioned complete-field
prediction on the observed finite reachable set.

The result also clarifies the contrast with the natural-feature
baselines. Among the 511 tested feature combinations, none achieved both
genuine field compression and exact prediction. This was not because
exact prediction intrinsically required nearly complete microscopic
field identification. Once the restriction to preselected features is
removed and the partition is derived from future incompatibilities, the
2,562 fields can be reduced to 692 classes. Hence,

\[{\boxed{\text{human-designed features}\neq\text{minimum distinctions required for prediction}}.}\]

The number 692 is not a universal state count for ABCW. It is the exact
solution under the following conditions:

\begin{itemize}
\tightlist
\item
  five agents;
\item
  the finite reachable set generated in this study;
\item
  the 56,536 observed transitions;
\item
  field-only compression;
\item
  retention of the action \({a}\);
\item
  one-step prediction; and
\item
  the complete \({\Delta W_{t+1}}\) as the prediction target.
\end{itemize}

Different initial conditions, reachable sets, numbers of agents,
prediction horizons, or coarser targets \({Q}\) may produce different
minimum partitions.

The theorem determines the chromatic number, and hence the minimum class
count of 692. The concrete map \({B_{692}}\) verified above is one
representative minimum partition obtained from the proper coloring
constructed by DSATUR. A different 692-coloring---including one
constructed from the anchor decomposition---may assign fields to classes
differently. We therefore claim uniqueness only of the minimum class
count for this incompatibility graph, not uniqueness of the minimizing
map or partition.

The search for the minimum number of classes is now complete.
Determining that the reduction is \({2562\to692}\), however, is not the
same as understanding what the 692 classes represent. Section 7
therefore examines which fields are merged within the classes, which
distinctions are exposed by action conditions, how the additional 67
distinctions are localized, and what symmetries are present in the
partition structure.

\begin{center}\rule{0.5\linewidth}{0.5pt}\end{center}

\hypertarget{structure-and-interpretation-of-the-minimum-partition}{%
\section{Structure and Interpretation of the Minimum
Partition}\label{structure-and-interpretation-of-the-minimum-partition}}

The preceding analysis established that 692 is the minimum number of
classes required to predict the next field \({\Delta W'}\) uniquely from
the current field \({\Delta W}\) conditioned on the action \({a}\), over
the 2,562 observed fields.

This result is not merely a statement that 2,562 states can be
compressed into 692 states. We now examine the structure of the 692
classes, why 625 classes are insufficient, and why 67 additional
distinctions are required.

The objective is not to assign an immediately interpretable meaning to
every class. Rather, it is to clarify how the distinctions required to
preserve the future relate to distinctions among the original
microscopic fields.

\hypertarget{the-minimum-as-a-predictive-distinction-count}{%
\subsection{The minimum as a predictive distinction
count}\label{the-minimum-as-a-predictive-distinction-count}}

The partition obtained here is not an arbitrary clustering of the
observed fields \({\Delta W}\). Two fields \({\Delta W_i}\) and
\({\Delta W_j}\) can be assigned to the same class only if identifying
them preserves unique prediction of the next field under all observed
action conditions shared by the two fields.

Among maps of the form

\[{B:\Delta W\longrightarrow\lbrace1,\ldots,K\rbrace,}\]

we require

\[{\bigl(a_t,B(\Delta W_t)\bigr)\longrightarrow\Delta W_{t+1}}\]

to be single-valued on the observed dataset. The exact minimum is

\[{K_{\min}=692.}\]

No universal meaning should be assigned to the number 692 itself. It is
determined by the five-agent ABCW model, the four initial fields, the
finite reachable set generated from them, and the target of predicting
the complete one-step-ahead field \({\Delta W_{t+1}}\). We therefore
refer to the result as the minimum predictive partition preserving
one-step, action-conditioned field responses on this dataset, rather
than as a universal number of ABCW macrostates.

\hypertarget{class-size-structure-and-genuine-merging}{%
\subsection{Class-size structure and genuine
merging}\label{class-size-structure-and-genuine-merging}}

The original observation set contains 2,562 distinct values of
\({\Delta W}\). The identity observation

\[{P_{\mathrm{id}}(\Delta W)=\Delta W}\]

distinguishes all of them and therefore has 2,562 classes. By contrast,
the minimum predictive partition reduces the distinction count to

\[{2562\longrightarrow692.}\]

In terms of field-class count, approximately 73\% of the original
distinctions can be discarded without losing one-step field prediction
on the observed dataset.

This compression is not produced solely by a small number of
exceptionally large classes. Of the 692 classes:

\begin{itemize}
\tightlist
\item
  352 contain a single \({\Delta W}\);
\item
  340 contain multiple distinct values of \({\Delta W}\); and
\item
  the largest class contains 65 distinct values of \({\Delta W}\).
\end{itemize}

Many fields must therefore remain individually distinguished, but
genuine merging is also widespread. Microscopically different fields can
often be treated as identical for the specified prediction task. The
result provides a concrete finite example in which

\[{\text{microscopically distinct}}\]

does not imply

\[{\text{must remain distinct for prediction}.}\]

\hypertarget{the-625-class-anchor-decomposition}{%
\subsection{The 625-class anchor
decomposition}\label{the-625-class-anchor-decomposition}}

The intermediate lower bound of 625 is central to the internal structure
of the exact result. It is obtained by fixing the anchor action

\[{a=(-,-,-,-,-).}\]

Under this action, 1,239 current fields are observed and produce 625
distinct next fields. Grouping current fields by their next field under
the anchor action gives an initial decomposition into 625 anchor future
groups.

This 625-group decomposition is not yet a valid complete predictive
partition when the other observed action conditions are included. Of the
625 groups, 575 require no further subdivision. All remaining conflicts
are localized in the other 50 groups:

\[{625=575+50.}\]

The increase from 625 to 692 is therefore not distributed uniformly
across the partition.

\hypertarget{the-additional-67-distinctions}{%
\subsection{The additional 67
distinctions}\label{the-additional-67-distinctions}}

Refining the 50 difficult groups so that future fields remain unique
under every observed action condition requires 67 additional
distinctions:

\[{625+67=692.}\]

After refinement, the 50 original groups split into

\[{39\times2+8\times3+2\times4+1\times7=117}\]

subclasses. The increase is therefore

\[{117-50=67.}\]

Most predictive distinctions are already forced at the 625-group stage,
whereas the remaining ambiguity is concentrated in a comparatively small
subset. The minimum partition thus has the form

\[{\text{large base partition}+\text{localized additional refinement}.}\]

This structure is not apparent if the partition is viewed only as a
lookup table containing 692 labels. Figure 6 summarizes both the final
composition of 352 singleton and 340 multi-field classes and the
localization of all 67 additional distinctions within 50 of the 625
anchor groups.

\begin{figure}
\centering
\includegraphics[trim=0 0 0 24bp,clip]{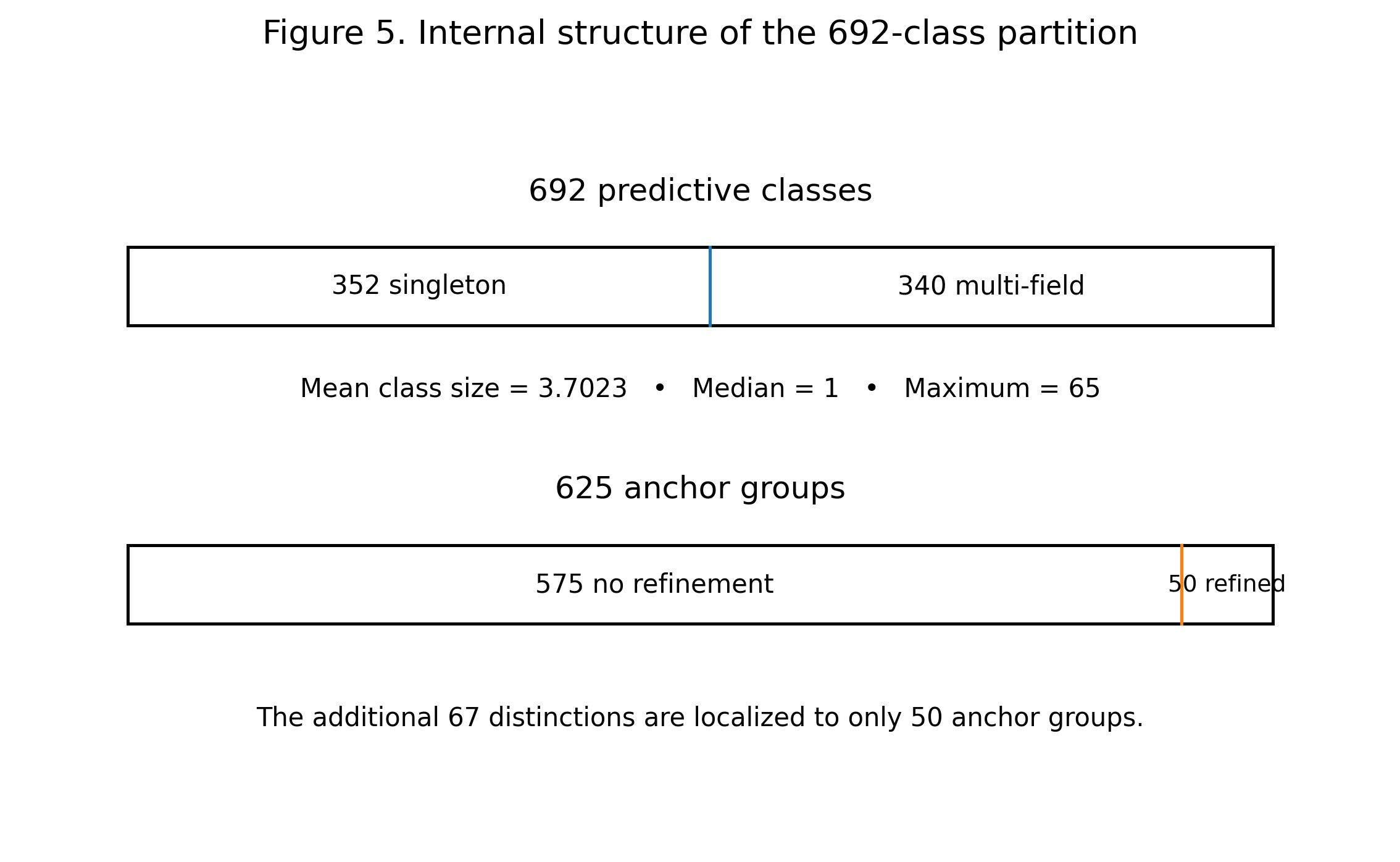}
\caption{Internal structure of the 692-class partition. Top:
composition of the 692 predictive classes (352 singleton and 340
multi-field classes), with verified mean, median, and maximum class
size. Bottom: localization of all additional refinements to 50 of the
625 anchor groups.}
\end{figure}

\hypertarget{geometric-proximity-versus-predictive-equivalence}{%
\subsection{Geometric proximity versus predictive
equivalence}\label{geometric-proximity-versus-predictive-equivalence}}

We next ask whether fields assigned to the same class are also similar
in an ordinary geometric sense. If the minimum predictive partition
resembled a simple geometric clustering, nearby values of \({\Delta W}\)
would be expected to share classes, while distant fields would tend to
be separated.

Inspection of the partition does not support this expectation in
general. Fields that differ substantially in their microscopic
representations can belong to the same predictive class, while nearby
fields can be assigned to different classes.

The partition therefore does not preserve the geometry of the
\({\Delta W}\) space itself. It preserves the action-conditioned
future-response structure

\[{\Delta W_t\longmapsto\left(a\longmapsto\Delta W_{t+1}\right).}\]

In this predictive sense, two fields are ``close'' not because their
matrix entries differ only slightly, but because they produce
indistinguishable future responses under the relevant observed action
conditions.

\hypertarget{predictive-sufficiency-and-representation-cost}{%
\subsection{Predictive sufficiency and representation
cost}\label{predictive-sufficiency-and-representation-cost}}

These results make concrete the relationship between predictive
sufficiency and compression. The identity observation
\({P_{\mathrm{id}}}\) necessarily retains the information required for
prediction, but it distinguishes all 2,562 fields and performs no
compression. At the opposite extreme, a very coarse observation may
compress the state description substantially while failing to
distinguish the relevant futures.

The question is therefore:

\begin{quote}
\textbf{How many distinctions can be discarded without losing the
specified future?}
\end{quote}

For the finite dataset and one-step field-prediction target considered
here, the boundary occurs at 692 classes. Compression and predictive
sufficiency must be treated as separate axes: coarseness alone is
insufficient, while predictability alone is trivially achieved by the
identity observation. The objective is to retain only the distinctions
required for prediction and discard the rest.

\hypertarget{why-the-partition-is-not-yet-an-interpretable-macrovariable}{%
\subsection{Why the partition is not yet an interpretable
macrovariable}\label{why-the-partition-is-not-yet-an-interpretable-macrovariable}}

The existence of a 692-class minimum partition does not establish that a
natural macrovariable of ABCW has been discovered.

First, the partition is constructed from an observed finite reachable
set. Second, its target is restricted to the one-step-ahead field
\({\Delta W_{t+1}}\). Third, the classes are defined by future
responses, and it has not been shown that they can be expressed through
simple physical or network quantities. Fourth, exact determination of
the minimum class count does not imply that the class labels or the
partition realizing that count are unique: the DSATUR-derived
\({B_{692}}\) is a verified representative optimum, not a canonical
quotient.

The result directly supports only the following claim:

\begin{quote}
\textbf{In the observed finite ABCW system, not all 2,562 microscopic
field distinctions are required to preserve exact one-step,
action-conditioned field responses; 692 predictive distinctions are
necessary and sufficient.}
\end{quote}

\hypertarget{an-intermediate-descriptive-level}{%
\subsection{An intermediate descriptive
level}\label{an-intermediate-descriptive-level}}

The 692 classes occupy an intermediate descriptive level between a
complete microscopic representation and a small set of simple
macroscopic statistics. They retain the distinctions required for a
specified future while discarding distinctions that are irrelevant to
that target on the observed dataset.

They should not, however, be treated as natural macroscopic variables of
ABCW in general. Their stability under larger systems, different initial
conditions, and longer prediction horizons remains unresolved, as does
the possibility of representing the 692 classes through a concise set of
network observables. Section 8 discusses these implications and
limitations.

\begin{center}\rule{0.5\linewidth}{0.5pt}\end{center}

\hypertarget{discussion}{%
\section{Discussion}\label{discussion}}

\hypertarget{what-the-exact-minimum-means}{%
\subsection{What the exact minimum
means}\label{what-the-exact-minimum-means}}

The central question of this study was:

\begin{quote}
\textbf{Which distinctions in the present must be retained in order to
predict the future?}
\end{quote}

The finite reachable set observed for the five-agent ABCW model contains
2,562 distinct current fields \({\Delta W}\). An identity field
representation that distinguishes all of them necessarily retains the
information required to predict the complete one-step-ahead field
\({\Delta W_{t+1}}\). The relevant question, however, is not whether the
complete current field is sufficient, but how many of its distinctions
can be discarded without losing the specified future.

The exact result obtained in Section 6 is

\[{|\operatorname{Im}B_{\min}|=\chi(G)=692.}\]

Under the conditions of the observed finite reachable set, field-only
compression, retention of the current action \({a_t}\), and exact
prediction of the complete next field \({\Delta W_{t+1}}\), the
current-field distinction count can be reduced as

\[{2562\longrightarrow692.}\]

With 691 or fewer classes, at least one pair of current fields producing
different next fields under a common action condition must be
identified, so exact prediction cannot be preserved.

The value 692 is therefore neither an arbitrary output of a compression
algorithm nor a performance score for a chosen feature set. It is the
minimum number of field distinctions that must be retained for the
specified prediction task.

This interpretation has explicit limits. The value is not the state
count of ABCW in general, nor is it a universal number of naturally
occurring macrostates. It depends on the number of agents, initial
fields, reachable set, observation form, prediction target, time
horizon, and requirement of exact prediction. The main result is
therefore most precisely stated as follows:

\begin{quote}
\textbf{In the finite ABCW system examined here, all 2,562 microscopic
field distinctions are unnecessary for predicting the complete
one-step-ahead field, but at least 692 predictive field distinctions
must be retained.}
\end{quote}

\hypertarget{microscopic-difference-and-predictive-relevance}{%
\subsection{Microscopic difference and predictive
relevance}\label{microscopic-difference-and-predictive-relevance}}

The central interpretation is that

\[{\text{microscopically different}}\]

and

\[{\text{must remain distinct for the specified prediction}}\]

are not equivalent.

The minimum partition contains 352 singleton classes and 340 classes
that merge multiple distinct fields, with the largest class containing
65 current fields. Many fields must remain individually identifiable,
but many others can be merged despite being microscopically different.

Class membership also does not coincide with simple geometric proximity
between field matrices. Fields that differ substantially can share a
predictive class, whereas nearby fields must be separated whenever they
produce different next fields under a common action condition.

The minimum partition therefore preserves not the static geometry of the
current fields themselves, but their observed action-conditioned
one-step response structure:

\[{\Delta W_t\longmapsto\left(a\longmapsto\Delta W_{t+1}\right).}\]

This demonstrates a limitation of evaluating the fineness of a state
description solely through the geometry or dimensionality of the present
state. Once a predictive objective is fixed, what matters is whether a
present distinction appears as a distinction in the target future.

\hypertarget{hand-designed-observables-and-dynamics-derived-partitions}{%
\subsection{Hand-designed observables and dynamics-derived
partitions}\label{hand-designed-observables-and-dynamics-derived-partitions}}

Before deriving the minimum partition, we evaluated field compression
based on human-selected network features, including norms, local
structure, and out-strength. These features are not uninformative. In
particular, out-strength gives a unique next field for 11,084 of 11,142
observed \({(a,\mathrm{out\mbox{-}strength})}\) states, corresponding to
a state-level determinism score of 99.4794\%.

Nevertheless, 58 conflicting observed states remain. Across all 511
nonempty combinations of the nine candidate feature families, none both
merged genuinely distinct values of \({\Delta W}\) and predicted the
complete next field exactly.

Viewed in isolation, this failure could suggest that exact prediction
requires an almost complete representation of the current field. Section
6 shows otherwise. Once the prior restriction to selected features is
removed and incompatibilities in future responses are used directly, the
2,562 fields can be compressed into 692 classes.

For this dataset,

\[{\text{human-designed features}\neq\text{minimum distinctions required for prediction}.}\]

This does not support a general conclusion that natural features are
unsuitable for prediction. The tested family is finite, and other
features or constructions may express the 692-class partition, or
another high-performing compression, concisely. The result instead
demonstrates a complementary direction of analysis: first derive from
the dynamics the distinctions required to preserve the future, and then
seek an interpretable representation of those distinctions.

The 692-class partition is therefore better understood as a predictive
benchmark for future searches for interpretable macroscopic descriptions
than as a final set of interpretable macrovariables.

\hypertarget{relation-to-state-aggregation-and-predictive-state-ideas}{%
\subsection{Relation to state aggregation and predictive-state
ideas}\label{relation-to-state-aggregation-and-predictive-state-ideas}}

As discussed in Section 2, the problem has connections to state
aggregation, lumpability, bisimulation, computational mechanics, and
minimization of incompletely specified finite-state machines.

The connection to lumpability lies in the shared objective of retaining
a well-defined description of future behavior after microscopic states
have been aggregated. In the present study, however, we distinguish
self-closure,

\[{P(X_t)\longrightarrow P(X_{t+1}),}\]

from prediction of a specified target,

\[{P(X_t)\longrightarrow Z.}\]

The preserved target is the complete one-step-ahead field
\({\Delta W_{t+1}}\), not the next state of the reduced field class. The
692-class partition is therefore not shown to define closed reduced
dynamics.

The study also shares with bisimulation and MDP minimization the
principle of aggregating states according to future behavior. It does
not, however, preserve an optimal policy, reward structure, or complete
reduced model. The current action is retained, only the field variable
is compressed, and the preserved output is the complete next field.

The conceptual connection to causal states in computational mechanics is
especially close because both approaches retain distinctions according
to their relevance for the future. Causal states classify histories by
conditional distributions over future sequences, whereas the present
construction partitions current fields in a finite reachable set and
restricts the target to one time step. Moreover, the input-output data
are partially specified: not every action condition is observed for
every current field. We therefore do not identify the 692 classes with
causal states or an \({\epsilon}\)-machine.

The most direct mathematical correspondence is with minimization of
incompletely specified finite-state machines. The data can be viewed as
the partially specified input-output relation

\[{(\Delta W,a)\longmapsto\Delta W'.}\]

Two current fields cannot share a class when they produce different
outputs under the same observed input condition. This
compatibility/incompatibility structure directly parallels the classical
problem. The present formulation does not require recursive consistency
of successor classes, because the preserved target is the complete
one-step-ahead field rather than a successor class label. Under this
depth-1 restriction, the minimum field-class count is given directly by
the chromatic number of the incompatibility graph.

The study therefore does not propose a new general theory of state
minimization or predictive representation. Its specific contribution is
to formulate this depth-1 predictive state-reduction problem for
partially specified data generated by a finite competitive agent system
and to solve the minimum exactly through a constructive upper bound and
an independent matching lower bound.

\hypertarget{the-structure-behind-692}{%
\subsection{The structure behind 692}\label{the-structure-behind-692}}

The exact minimum is not produced by uniformly refining all 2,562
fields. Under the anchor action

\[{a=(-,-,-,-,-),}\]

the 1,239 observed current fields produce 625 distinct next fields,
forcing an initial 625-group decomposition. When the other observed
action conditions are considered, 575 groups contain no additional
internal conflicts, while 50 require further refinement.

Exact refinement of these 50 groups produces 117 subclasses, increasing
the total class count by

\[{117-50=67.}\]

Thus,

\[{625+67=692.}\]

The predictive distinctions are not distributed uniformly across current
fields. For most anchor groups, the distinction exposed by a single
action condition is sufficient. Additional action conditions require new
distinctions only within 50 groups. The minimum partition can therefore
be understood as

\[{\text{large base partition}+\text{localized additional refinement}.}\]

This localization does not establish a general law that predictive
information is localized. It is an observed structural property of the
present finite dataset and prediction task. Nevertheless, it shows that
tracing which conditions expose new incompatibilities can reveal where
predictively relevant differences occur, information that is hidden if
the result is represented only as a lookup table of 692 labels.

\hypertarget{scope-and-limitations}{%
\subsection{Scope and limitations}\label{scope-and-limitations}}

The result has at least six important limitations.

First, the analysis concerns a five-agent ABCW model and the finite
reachable set generated from four initial fields. The value 692 cannot
be extrapolated directly to other numbers of agents, other initial
conditions, or ABCW systems in general.

Second, the current action \({a_t}\) is retained, and only the field
variable \({\Delta W_t}\) is compressed. The study does not determine a
general minimum state representation that compresses action and field
jointly.

Third, the prediction target is the complete one-step-ahead field
\({\Delta W_{t+1}}\). We have not established that the 692 classes are
sufficient for multistep prediction or that they support closed
recursive dynamics.

Fourth, the analyzed input-output behavior is partially specified. Not
every action condition is observed for every current field. The exact
minimum therefore applies to the observed finite reachable set and does
not guarantee predictions for unobserved combinations.

Fifth, exact one-step prediction is imposed as a hard constraint. This
makes the minimum partition mathematically well defined, but it does not
imply that exact prediction is the only scientifically meaningful
objective. Simpler representations such as out-strength can provide
high, though imperfect, predictive performance.

Sixth, establishing a minimum of 692 classes does not show that those
classes can be expressed through a concise set of network features, an
algebraic rule, or a small number of interpretable variables. Nor does
the minimum class count imply uniqueness of the partition that attains
it.

These limitations specify, rather than weaken, the scope of the main
result. What has been established exactly is the minimum number of field
classes necessary and sufficient under a specified finite dataset,
observation form, prediction target, and time horizon.

\hypertarget{future-directions}{%
\subsection{Future directions}\label{future-directions}}

Several extensions follow naturally from these limitations.

The first is to extend the prediction horizon. Two fields that can be
merged for one-step prediction may generate different futures over
multiple time steps. Conversely, if the target is changed from the
complete future field to a coarser quantity, a partition with fewer than
692 classes may suffice.

The observation class may also be generalized from field-only
compression to a joint map

\[{C:(a,\Delta W)\longrightarrow\mathcal Z,}\]

so that action distinctions and field distinctions are optimized
simultaneously rather than retaining \({a}\) as an uncompressed
conditioning variable. For an \({h}\)-step target, depth-1 compatibility
is no longer sufficient: on the partially specified transition system,
two records merged at depth \({h}\) must agree on the required immediate
output and, whenever both successors are specified, their successor
records must remain compatible at depth \({h-1}\). This gives the
recursive compatibility condition familiar from incompletely specified
FSM minimization and provides a direct route to finite-horizon
extensions of the present graph construction.

The second is to move from exact to approximate prediction. This study
imposes exact prediction and minimizes the resulting class count. In
practical model analysis, the trade-offs among predictive performance,
compression, and representation cost may be more informative. Many
intermediate representations may exist between simple, highly predictive
features and the exact 692-class partition. These three quantities
should therefore be treated as separate axes:

\[{\text{predictive performance},}\]

\[{\text{compression},}\]

\[{\text{representation cost}.}\]

Studying their Pareto structure is a natural next step.

The third direction is to seek a readable representation of the
692-class partition. The present result shows that the distinctions
required for prediction can be derived from the dynamics, but it remains
unknown whether they can be re-expressed through a small set of
interpretable network quantities.

The fourth is to vary the number of agents, the initial fields, and the
reachable set. Such comparisons are required to determine which aspects
of the observed structure reflect more general properties of ABCW and
which are specific to the present finite dataset.

These extensions are outside the main result of this paper. What is
solved here is the more limited depth-1 problem. That restriction makes
it possible to turn the abstract question of which distinctions matter
for prediction into a finite problem that can be inspected and solved
exactly.

\begin{center}\rule{0.5\linewidth}{0.5pt}\end{center}

\hypertarget{conclusion}{%
\section{Conclusion}\label{conclusion}}

\hypertarget{main-result-1}{%
\subsection{Main result}\label{main-result-1}}

This study determined the minimum number of field classes required, on a
finite reachable set of the five-agent ABCW model, to preserve exact
prediction of the complete one-step-ahead field \({\Delta W_{t+1}}\)
while retaining the current action \({a_t}\) and compressing only the
current field \({\Delta W_t}\).

We constructed an incompatibility graph on the 2,562 observed current
fields. A proper 692-coloring provides a constructive upper bound, while
an independent lower bound derived from an anchor action requires the
same number of colors. Hence,

\[{|\operatorname{Im}B_{\min}|=\chi(G)=692.}\]

Not all 2,562 microscopic field distinctions are therefore required for
the specified prediction task. At the same time, no field-only partition
with 691 or fewer classes can preserve exact prediction.

\hypertarget{interpretation-and-scope}{%
\subsection{Interpretation and scope}\label{interpretation-and-scope}}

The exact minimum is not produced by a uniform refinement of the field
set. The anchor action forces 625 base groups, of which only 50 require
further subdivision under other action conditions. These localized
refinements add 67 classes,

\[{625+67=692,}\]

and the final partition contains 352 singleton classes and 340
multi-field classes.

The natural-feature baselines also clarify the significance of the
result. Across 511 combinations of nine feature families, the
out-strength vector achieved a state-level determinism score of
99.4794\%, but no tested candidate combined genuine field compression
with exact prediction. In this dataset, the features that appeared
natural in advance did not coincide with the minimum distinctions
required to preserve the future.

The value 692 is not a universal number of macrostates for ABCW. It
depends on the five-agent model, the finite reachable set examined here,
retention of the action, field-only compression, and exact one-step
field prediction. Determining the minimum class count also does not
imply uniqueness of the minimum coloring or establish that the classes
can be expressed through a concise analytic observation map.

\hypertarget{outlook}{%
\subsection{Outlook}\label{outlook}}

Natural extensions include multistep prediction, trade-offs between
compression and approximate predictive performance, scaling with the
number of agents and the choice of initial fields, and the search for a
smaller set of interpretable structural quantities capable of
representing the 692-class partition.

The most limited conclusion of this paper is that \textbf{being distinct
as a state is not the same as needing to remain distinct for a specified
prediction}. The 692-class minimum predictive partition of this finite
ABCW system provides one exact answer to the question of how much of the
present can be forgotten without losing the specified future.

\begin{center}\rule{0.5\linewidth}{0.5pt}\end{center}

\hypertarget{data-and-code-availability}{%
\section{Data and Code Availability}\label{data-and-code-availability}}

The complete transition dataset, source code, exact-coloring
certificate, lower-bound reconstruction, natural-feature analysis, and
automated verification scripts are available at
\url{https://github.com/youque8/abcw-paper-reproducibility}. The
repository release corresponding to this manuscript is version 1.1.0.
The core numerical chain can be checked by running
\texttt{python\ scripts/verify\_all.py}; the full recomputation,
including all 511 natural-feature combinations and exact re-solving of
the anchor refinements, is available through
\texttt{python\ scripts/verify\_full.py}.

\hypertarget{references}{%
\section*{References}\label{references}}
\addcontentsline{toc}{section}{References}

\hypertarget{refs}{}
\begin{CSLReferences}{1}{0}
\leavevmode\vadjust pre{\hypertarget{ref-arthur1994}{}}%
Arthur, W. Brian. 1994. {``Inductive Reasoning and Bounded
Rationality.''} \emph{American Economic Review} 84 (2): 406--11.

\leavevmode\vadjust pre{\hypertarget{ref-brelaz1979}{}}%
Brélaz, Daniel. 1979. {``New Methods to Color the Vertices of a
Graph.''} \emph{Communications of the ACM} 22 (4): 251--56.
\url{https://doi.org/10.1145/359094.359101}.

\leavevmode\vadjust pre{\hypertarget{ref-bron1973}{}}%
Bron, Coen, and Joep Kerbosch. 1973. {``Algorithm 457: Finding All
Cliques of an Undirected Graph.''} \emph{Communications of the ACM} 16
(9): 575--77. \url{https://doi.org/10.1145/362342.362367}.

\leavevmode\vadjust pre{\hypertarget{ref-challet1997}{}}%
Challet, Damien, and Yi-Cheng Zhang. 1997. {``Emergence of Cooperation
and Organization in an Evolutionary Game.''} \emph{Physica A:
Statistical Mechanics and Its Applications} 246 (3--4): 407--18.
\url{https://doi.org/10.1016/S0378-4371(97)00419-6}.

\leavevmode\vadjust pre{\hypertarget{ref-givan2003}{}}%
Givan, Robert, Thomas Dean, and Matthew Greig. 2003. {``Equivalence
Notions and Model Minimization in {M}arkov Decision Processes.''}
\emph{Artificial Intelligence} 147 (1--2): 163--223.
\url{https://doi.org/10.1016/S0004-3702(02)00376-4}.

\leavevmode\vadjust pre{\hypertarget{ref-kemeny1960}{}}%
Kemeny, John G., and J. Laurie Snell. 1960. \emph{Finite Markov Chains}.
Van Nostrand.

\leavevmode\vadjust pre{\hypertarget{ref-kohavi2009}{}}%
Kohavi, Zvi, and Niraj K. Jha. 2009. \emph{Switching and Finite Automata
Theory}. 3rd ed. Cambridge University Press.

\leavevmode\vadjust pre{\hypertarget{ref-paull1959}{}}%
Paull, Marvin C., and Stephen H. Unger. 1959. {``Minimizing the Number
of States in Incompletely Specified Sequential Switching Functions.''}
\emph{IRE Transactions on Electronic Computers} EC-8 (3): 356--67.
\url{https://doi.org/10.1109/TEC.1959.5222697}.

\leavevmode\vadjust pre{\hypertarget{ref-shalizi2001}{}}%
Shalizi, Cosma Rohilla, and James P. Crutchfield. 2001. {``Computational
Mechanics: Pattern and Prediction, Structure and Simplicity.''}
\emph{Journal of Statistical Physics} 104 (3--4): 817--79.
\url{https://doi.org/10.1023/A:1010388907793}.

\leavevmode\vadjust pre{\hypertarget{ref-tishby2000}{}}%
Tishby, Naftali, Fernando C. Pereira, and William Bialek. 2000. {``The
Information Bottleneck Method.''}
\url{https://arxiv.org/abs/physics/0004057}.

\end{CSLReferences}

\end{document}